\documentclass{siamltex}
\usepackage{amsmath,amssymb,amsfonts,latexsym,graphicx,color,multirow,footnote}

\newtheorem{algorithm}{Algorithm}

\usepackage{amsmath,amsfonts,amssymb}
\usepackage[dvipsone]{epsfig}
\usepackage{ulem}

\title{A FEAST variant
incorporated
with a power
iteration
 }

\author{
Man-Chung Yeung\thanks{Department of Mathematics,
University of Wyoming, Laramie, WY, USA ({\tt myeung, llee@uwyo.edu}).}
\and Long Lee$^*$
}

\begin{document}
\maketitle

\begin{abstract}
We present a variant of the FEAST matrix eigensolver for solving restricted real and symmetric
eigenvalue problems. The method is derived from a combination of
the FEAST method
and a power subspace iteration
process. Compared with the original FEAST method, our
new method does not require that the search subspace dimension must be greater than or equal to the number of eigenvalues in
a search interval. Together with two contour integrations per iteration, the new method can deal with relatively narrow search intervals more effectively.
Empirically, the FEAST iteration and the power subspace iteration are in a mutually beneficial collaboration to make
the new method stable and robust.
\end{abstract}

\begin{keywords}
FEAST eigensolver, power iteration,
contour integral, spectral projection.
\end{keywords}

\begin{AMS}
15A18, 58C40, 65F15
\end{AMS}

\section{Introduction}\label{sec:intro}
Consider the eigenvalue problem
\begin{equation}\label{equ:7-14-1}
\begin{array}{ccc}
A x = \lambda x, & & x \ne 0,
\end{array}
\end{equation}
where $A \in {\mathbb R}^{n \times n}$ is symmetric. Since $A$ is real and symmetric, its eigenvalues $\lambda$ are real numbers. Given an open interval $(a, b)$ on the real axis, our aim is to extract the first $l$
largest eigenvalues from the set
$
\{\lambda | \lambda \in (a, b)\}$,
along with their associated eigenvectors.

We solve the problem
by FEAST\cite{polizzi} combined with a conventional power subspace iteration process.
Let ${\cal K}$ denote the eigenspace of $A$ associated with
the eigenvalues in $(a, b)$, and pick a random
$Y \in {\mathbb R}^{n \times m}$ such that
\begin{equation}\label{equ:10-16-1}
span\{Y\} \subset {\cal K},
\end{equation}
where $span\{Y\}$ is the column space of $Y$ and $m \geq l$.
The power subspace iteration, when applied to $A - \sigma I$ for some appropriately chosen shift $\sigma$ and starting with $Y$, will produce approximations to the
first $l$ largest eigenvalues
of $A$ in $(a, b)$ and their associated eigenvectors by projecting the problem onto the column space of $(A - \sigma I)^k Y$ in its $k$th iteration.
Theoretically we should have
\begin{equation}\label{equ:10-16-2}
span\{(A - \sigma I)^k Y\} \subset {\cal K}
\end{equation}
for all $k$. Computationally, however, property (\ref{equ:10-16-2}) is gradually lost as $k$ is increasing, and some correction is needed from time to time to keep
(\ref{equ:10-16-2}) hold as much as possible.

Suppose $Y_k$ is a basis matrix of the subspace $span\{(A - \sigma I)^k Y\}$. One standard correction on $Y_k$ is to compute a Cauchy integral of the form
\begin{equation}\label{equ:10-17-1}
Z_k = \frac{1}{2 \pi i} \oint_\Gamma (zI - A)^{-1} Y_k dz,
\end{equation}
where $i = \sqrt{-1}$ and $\Gamma$ is the counterclockwise oriented circle in the complex plane with its center at $c = (a + b)/2$ and radius $r = (b-a)/2$.
The computation of the integral
requires the solution of a bunch of linear systems arising from
discretizing the integral (\ref{equ:10-17-1}) by a quadrature rule.
When $(a, b)$ is narrow relative to the spectrum of $A$, which results in a relatively small $r$, the linear systems can be hard to solve
because the poles of the resulting rational filter are too close to the real axis.
To overcome the problem, instead of using (\ref{equ:10-17-1}) we adopt the strategy for using integrals (\ref{equ:10-18-1}) as a corrector for $Y_k$:
\begin{equation}\label{equ:10-18-1}
\begin{array}{ccc}
\displaystyle{W_k = \frac{1}{2 \pi i} \oint_{\Gamma_L} (zI - A)^{-1} Y_k dz}, & & \displaystyle{Z_k = \frac{1}{2 \pi i} \oint_{\Gamma_R} (zI - A)^{-1} W_k dz},
\end{array}
\end{equation}
where $\Gamma_L$ and $\Gamma_R$ are two counterclockwise oriented circles in the complex plane that have equal radii of $r$ with $r \geq b-a$. The center of
the left circle $\Gamma_L$ is
$c_L = b - r$ and that of the right one $\Gamma_R$ is $c_R = a + r$.
The overlap of the circles on the real axis is exactly the interval $(a, b)$. Theoretically we have $span\{Z_k\} \subset {\cal K}$ for the $Z_k$ in
(\ref{equ:10-18-1}) even though $span\{Y_k\} \not\subset {\cal K}$. The advantage of (\ref{equ:10-18-1}) over (\ref{equ:10-17-1}) is that the radius $r$
can be arbitrarily chosen
provided that $r \geq b-a$.
As a result, we can pick a relatively large $r$ to avoid ill-conditioned linear systems to
occur.

The approach described above has led to an algorithm that can be viewed as a combination of the FEAST subspace iteration,
a spectral projection subspace iteration process,
and a conventional power subspace iteration.
We name the algorithm {\it FEAST-power subspace iteration with two contour integrations
per iteration}, or
abbreviated as F$_2$P.

The FEAST method was developed by Polizzi in \cite{polizzi}
to compute all the eigenvalues of (\ref{equ:7-14-1}) inside a given interval $(a, b)$, and their corresponding eigenvectors. A software package can be found in \cite{polizzi_software}.
FEAST is a subspace iteration method with a Rayleigh-Ritz procedure and a spectral projection
procedure.
Its stability and robustness have been demonstrated in \cite{kramer} and in other applications. Theoretical analysis exists in \cite{TP13} and a comparison
with some existing
Krylov subspace eigensolvers was made in \cite{GP}. The numerical computation and analysis of rational filters from (\ref{equ:10-17-1}) have been discussed in depth in \cite{GPTV, yxccb},
and filters other than those induced by (\ref{equ:10-17-1}) were proposed in \cite{xisaad}
through a least-squares process and in \cite{bib:WD2019} through an optimization
process.
Moreover,
generalizations from the symmetric or Hermitian case to the non-Hermitian or even generalized non-Hermitian case have been made in \cite{kpt, yxccb, yin2019, ycy}.

Compared with FEAST, F$_2$P has several advantages in computation:
(i) it handles a narrow interval that contains the wanted eigenvalues in a more effective way. Here the narrowness is relative to the
spectrum of the matrix; (ii) it has more freedom in the choice of circle radius $r$. In theory, $r$ can be any number with $r \geq b -a$;
(iii) it removes the restriction $m \ge s$ on FEAST where
$s$ the number of eigenvalues inside the interval of interest.
On the other hand, a clear disadvantage for F$_2$P
is that
the computational cost could be higher, compared to FEAST, due to the use of
two contour integrations per iteration. Nonetheless, the robustness of F$_2$P may be able to compensate this disadvantage.

To aid the reader, we now outline the contents of the remainder of this paper. In \S\ref{sec:feast},
we briefly review the FEAST algorithm by employing the general subspace iteration setting in \cite{TP13}.
In \S\ref{sec:F2P}, we develop
a F$_2$P algorithm.
In \S\ref{sec:numexp}, numerical experiments are reported to illustrate the
robustness and applicability of the F$_2$P algorithm. Finally,
conclusions are made in \S\ref{sec:conclusions}. We also note that throughout the paper,
algorithms are presented in \textsc{Matlab} style. \textsc{Matlab} functions are written in typewriter font.

\section{The FEAST method}\label{sec:feast} Consider the eigenvalue problem (\ref{equ:7-14-1}). Given an open interval $(a, b)$ on the real axis,
we want to compute all or some of
the eigenvalues inside $(a, b)$ together with their associated eigenvectors.
This restricted eigenproblem is solved through
a power subspace iteration with the Rayleigh-Ritz procedure described in \cite{TP13}. A slightly modified version of the power subspace iteration with $B = I$ is presented below.

\begin{algorithm}\label{alg:9-4-1} A general power subspace iteration with Rayleigh-Ritz
\end{algorithm}
\begin{tabbing}
x\=xxx\= xxx\=xxx\=xxx\=xxx\=xxx\=xxx\=xxx\kill
\>1. \> Pick $Y_{(0)} \in {\mathbb R}^{n \times m}$ randomly. Set $k \leftarrow 1$.\\
\>2.\> {\bf repeat}\\
\>3.\>\> $Q_{(k)} \leftarrow \rho(A) Y_{(k-1)}$.\\
\>4.\>\> $\hat{A}_{(k)} \leftarrow Q_{(k)}^T A Q_{(k)}$, $\hat{B}_{(k)} \leftarrow Q_{(k)}^T Q_{(k)}$.\\
\>5.\>\> Solve the $m$-dimension eigenproblem $\hat{A}_{(k)} \hat{X}_{(k)} = \hat{B}_{(k)} \hat{X}_{(k)} \hat{\Lambda}_{(k)}$ for $\hat{\Lambda}_{(k)}, \hat{X}_{(k)}$.\\
\>6.\>\> $Y_{(k)} \leftarrow Q_{(k)} \hat{X}_{(k)}$.\\
\>7.\>\> $k \leftarrow k + 1$.\\
\>8.\> {\bf until} Appropriate stopping criteria.
\end{tabbing}

\vspace{.2cm}

The $\rho$ in the algorithm is a mapping from
${\mathbb R}^{n \times n}$ to ${\mathbb R}^{n \times n}$.
In the case when $\rho(A) = A$, the above algorithm is the standard power subspace iteration
(see, for instance, Algorithm 5.3 in \cite{saad} and the algorithm in Table 14.2 in \cite{Parlett}). On the other hand, if we denote by $Q_{(a, b)}$ the orthonormal  matrix of eigenvectors associated with the eigenvalues of $A$
in $(a, b)$ and set $\rho(A) = Q_{(a, b)} Q_{(a, b)}^T$, then the algorithm is the FEAST algorithm. Further,
if we set $\rho(A) = Q_{(a, b)} (Q_{(a, b)}^T A Q_{(a, b)}) Q_{(a, b)}^T$,\footnote{It can be seen that $Q_{(a, b)} (Q_{(a, b)}^T A Q_{(a, b)}) Q_{(a, b)}^T =
A Q_{(a, b)} Q_{(a, b)}^T = Q_{(a, b)} Q_{(a, b)}^T A$.}
the algorithm is the F$_2$P algorithm whose implementation version is presented in Algorithm \ref{alg:7-29-1} in \S\ref{subsec: FFS}.

$Q_{(a, b)} Q_{(a, b)}^T$ is a spectral projector onto the invariant eigenspace ${\mathbb K}$. To obtain an approximation to this projector,
one usually constructs a rational filter related to the interval $(a, b)$ and apply the filter to $A$, e.g., the approximate projectors from Zolotarev rational filters \cite{GPTV, Viaud} and from least-squares rational
filters \cite{xisaad}. In this paper, we adopt the one, used and studied in \cite{polizzi, TP13} and briefly described below, from a filter obtained by applying Gauss-Legendre quadrature rule to a Cauchy integral.

Let $\Gamma$ be the positively oriented circle in the complex plane with center at $c = (a+b)/2$ and radius $r = (b-a)/2$. The residue
\begin{equation}\label{equ:9-26-1}
P_\Gamma = \frac{1}{2 \pi i} \oint_\Gamma (z I - A)^{-1}dz
\end{equation}
then defines a projection operator onto the eigenspace ${\mathbb K}$ (see, for instance, \cite{polizzi, saad, TP13}), where $i = \sqrt{-1}$.
In fact, it can be shown that
\begin{equation}\label{equ:7-20-1}
P_\Gamma = Q_{(a,b)} Q^T_{(a,b)}.
\end{equation}
The contour integral in (\ref{equ:9-26-1}) is usually evaluated approximately by using a quadrature rule. To the end, we define the change of variable
$$
\begin{array}{ccc}
z = c + r e^{i \pi t}, & & -1 \leq t < 1.
\end{array}
$$
Then (\ref{equ:9-26-1}) is transformed to
$$
\begin{array}{rl}
P_\Gamma 
&= \displaystyle{\frac{r}{2} \int_{-1}^{1} e^{i \pi t} [(c + r e^{i \pi t}) I - A]^{-1} \,dt}\\ \\
&= \displaystyle{\frac{r}{2}\left [\int_{-1}^{0} e^{i \pi t} [(c + r e^{i \pi t}) I - A]^{-1} \,dt + \int_{0}^{1} e^{i \pi t} [(c + r e^{i \pi t}) I - A]^{-1} \,dt\right]}\\ \\
&= \displaystyle{r \int_{0}^{1} \mbox{Real}\{e^{i \pi t} [(c + r e^{i \pi t}) I - A]^{-1}\} \,dt},
\end{array}
$$
where we have exploited the fact that $c, r$, and $A$ are real in the final equation.
%
The integral in
the final equation is now approximated by using for example the Gauss-Legendre quadrature rule \cite{DR84}
on the interval $[0, 1]$ with truncation order $q$:
$$
P_\Gamma \approx r\sum_{k = 1}^q \omega_k \, \mbox{Real}\{e^{i \pi t_k} [(c + r e^{i \pi t_k}) I - A]^{-1}\} \equiv \zeta(A),
$$
where $\omega_k$ and $t_k$ are the weights and nodes on $[0, 1]$.

In practice, the $\rho(A)$ in Algorithm \ref{alg:9-4-1} is replaced with $\zeta(A)$
in the case of the FEAST algorithm and with $A \zeta(A)$ in the case of the F$_2$P algorithm. The computation of $\zeta(A) Y_{(k-1)}$ is
the dominant cost in the two algorithms.

Now for any $Y \in {\cal R}^{n \times m}$,
\begin{equation}\label{equ:7-15-2}
P_\Gamma Y \approx
\zeta(A) Y = r\sum_{k = 1}^q \omega_k \, \mbox{Real}\{e^{i \pi t_k} [(c + r e^{i \pi t_k}) I - A]^{-1} Y\}.
\end{equation}
There are $mq$ linear systems to solve in (\ref{equ:7-15-2})
\begin{equation}\label{equ:7-15-3}
\begin{array}{ccc}
[(c + r e^{i \pi t_k}) I - A] x = y_j, & & k = 1, \ldots, q,\,\, j = 1,\ldots, m.
\end{array}
\end{equation}
The linear systems, however, are independent of each other and can be solved in parallel.

The condition numbers of the coefficient matrices of the linear systems in (\ref{equ:7-15-3})
depend on $r, t_k$, and $A$,
but not on $c$.
Suppose all the eigenvalues of $A$ are contained by the interval $[x_0 - \Delta, x_0 + \Delta]$ for some $x_0 \in {\mathbb R}$ and $\Delta > 0$, and let $z_k
= c + r e^{i \pi t_k} = c + r \cos(\pi t_k) + i r \sin(\pi t_k) \equiv \alpha_k + i\beta_k$.
  Assume that $\alpha_k \in [x_0 - \Delta, x_0 + \Delta]$. Let $M = \alpha_k I - A$. Then $z_kI - A = M + i\beta_k I$, and $\sigma(M) \subset
[(\alpha_k - x_0) -\Delta, (\alpha_k - x_0) +\Delta]$.
Since
$$
(z_kI - A)^H (z_kI - A) = (i \beta_k I + M)^H (i \beta_k I + M)
= \beta_k^2 I  + M^2,
$$
the largest and the smallest singular values of $z_k I - A$ satisfy
    $$
    \begin{array}{l}
    \sigma_{max}(z_k I - A) = [\lambda_{max}(\beta_k^2 I  + M^2)]^{1/2} = [\beta_k^2  + \lambda_{max}(M^2)]^{1/2}
    \leq [\beta_k^2  + 4\Delta^2]^{1/2} \\ \\
    \sigma_{min}(z_k I - A) = [\lambda_{min}(\beta_k^2 I  + M^2)]^{1/2} = [\beta_k^2  + \lambda_{min}(M^2)]^{1/2}
    \geq [\beta_k^2  + 0]^{1/2} = \beta_k
    \end{array}
    $$
    Thus the condition number of $z_k I - A$ can be bounded as follows
\begin{equation}\label{equ:7-18-1}
\begin{array}{rl}
\kappa_2 (z_k I - A) & \displaystyle{= \frac{\sigma_{max}(z_k I - A)}{\sigma_{min}(z_k I - A)} \leq \frac{(\beta_k^2  + 4\Delta^2)^{1/2}}{\beta_k} = \left(1  + 4\frac{\Delta^2}{\beta_k^2}\right)^{1/2}}\\ \\
& \displaystyle{\leq 1  + 2\frac{\Delta}{\beta_k}
= 1  + 2\frac{\Delta}{r \sin(\pi t_k)}}.
\end{array}
\end{equation}
This bound shows that, the larger the radius $r$ is and the farther away the $t_k$'s stay from the endpoints of the interval $[0, 1]$,
the better-conditioned the linear systems in (\ref{equ:7-15-3}) are.

The FEAST algorithm is numerically stable. It can catch the desired eigenvalues and eigenvectors accurately when it converges (see, for instance, \cite{kramer, polizzi}). The $\rho(A) Y_{(k-1)}$ in Algorithm \ref{alg:9-4-1} obtained through (\ref{equ:7-15-2}) is just an approximation, and there is a distance between ${\mathbb K}$ and
$span(\zeta(A) Y_{(k-1)})$. The distance, however, attenuates exponentially through the iteration process in the algorithm (see \cite{TP13} for the detail).

The computational cost of the FEAST algorithm is mainly in the solution of the linear systems in (\ref{equ:7-15-3}), where $A$ is usually large and sparse. As indicated in \cite{GP}, an optimized sparse direct solver (such as PARDISO \cite{SGKRH}) is typically used to solve the linear systems. Krylov subspace solvers, on the other hand, are also applicable and have been studied systematically in \cite{GP}.

We can observe two challenges about the implementation of the FEAST algorithm. First, if the provided interval $(a, b)$ in which the eigenpairs are desired is narrow
relative to the spectrum of $A$ (precisely, relative to $\Delta$), the
radius $r = (b - a)/2$ of the circle $\Gamma$ is small. In this case, the linear systems in (\ref{equ:7-15-3}) are likely to be ill-conditioned to solve according to (\ref{equ:7-18-1}). Of course, one can choose a larger interval $(\hat{a}, \hat{b})$ containing $(a, b)$ and compute the eigenvalues in $(\hat{a}, \hat{b})$, then extract those in $(a, b)$, but then some extra computational cost is required and the cost may not be modest.
Second, FEAST will fail to converge if the column size $m$ of the starting matrix $Y_{(0)}$ is less than the exact number $s$ of the eigenvalues in the interval $(a, b)$.
In other words, that $m \geq s$ is a necessary condition for FEAST to converge.
So $m$ depends on $s$ strongly.

Noting the challenges, in the next section, we propose solutions
to overcome the difficulties. Our solutions answer the following questions
(i) can we choose a large $r$ in the case when the interval $(a, b)$ is relatively small? (ii) can we release the restriction of $m \geq s$ from the FEAST algorithm?

\section{The FEAST-power subspace iteration method}\label{sec:F2P} The answers to the questions at the end of
\S\ref{sec:feast} lie in the following observations.
\begin{enumerate}
\item[(1)] Observation for question (i): The contour integral (\ref{equ:9-26-1}) on
 $(zI - A)^{-1}$ over a circle that encloses exactly the desired eigenvalues is a projection operator onto the associated eigenspace ${\mathbb K}$. However, a combination of contour integrals on $(zI - A)^{-1}$ over two circles whose overlapping region contains exactly the same desired eigenvalues will also provide a projection operator onto ${\mathbb K}$. Theoretically the two
circles can be chosen arbitrarily large. In \S \ref{sec:F_2}, we introduce FEAST$_2$, a FEAST algorithm with two contour integrations, that allows one to choose a large $r$.

\item[(2)] Observation for question (ii): Suppose the interval $(a, b)$ contains the dominant eigenvalues for$A$. As $k$ is increased and with some appropriate normalization
on $A^k$, the dominant eigenvalues of $A^k$ remain in $(a, b)$, but the relatively small eigenvalues are leaving the interval. As a result, the number of eigenvalues of $A^k$ in $(a, b)$ is decreasing as $k$ is increasing.
Thus we can apply the FEAST algorithm to $A^k$ for large enough $k$'s with a relatively small $m$. In \S \ref{subsec: subiter}, we explain the Power Subspace Iteration algorithm
and in \S \ref{subsec: FFS}, we combine
PSI and FEAST$_2$ to obtain a  FEAST$_2$-PSI algorithm.
\end{enumerate}

\subsection{FEAST with two contour integrations}\label{sec:F_2} Pick two circles $\Gamma_L$ and $\Gamma_R$ in the complex plane, as described in \S\ref{sec:intro}.
Define
$$
\begin{array}{ccc}
\displaystyle{P_{\Gamma_L} = \frac{1}{2 \pi i} \oint_{\Gamma_L} (z I - A)^{-1}dz} & \mbox{and} &
\displaystyle{P_{\Gamma_R} = \frac{1}{2 \pi i} \oint_{\Gamma_R} (z I - A)^{-1}dz}.
\end{array}
$$
According to (\ref{equ:7-20-1}), $P_{\Gamma_L} = Q_{(b - 2r, b)} Q_{(b - 2r, b)}^T$ and $P_{\Gamma_R} = Q_{(a, a+2r)} Q_{(a, a+2r)}^T$. With an appropriate rearrangement of the columns of $Q_{(b - 2r, b)}$ and $Q_{(a, a+2r)}$, we can express them as
$$\begin{array}{ccc}
Q_{(b - 2r, b)} = [Q_{(b - 2r, a]}, Q_{(a, b)}], & & Q_{(a, a+2r)} = [Q_{(a, b)}, Q_{[b, a+2r)}]
\end{array}
$$
Thus
$$
Q_{(a, a+2r)}^T Q_{(b - 2r, b)} = \left[ \begin{array}{c}
Q_{(a, b)}^T\\
Q_{[b, a+2r)}^T
\end{array} \right] [Q_{(b - 2r, a]}, Q_{(a, b)}] = \left[\begin{array}{cc}
0& I\\
0&0
\end{array}\right]
$$
and hence
$$
\begin{array}{rl}
P_{\Gamma_L \cap \Gamma_R} \equiv P_{\Gamma_R} P_{\Gamma_L} &= (Q_{(a, a+2r)} Q_{(a, a+2r)}^T) (Q_{(b - 2r, b)} Q_{(b - 2r, b)}^T)\\
& = Q_{(a, a+2r)} \left[\begin{array}{cc}
0& I\\
0&0
\end{array}\right] Q_{(b - 2r, b)}^T\\
& = Q_{(a, b)} Q_{(a, b)}^T.
\end{array}
$$
which is the orthogonal projector in (\ref{equ:7-20-1}).

The following algorithm is an implementation version of FEAST running with $P_{\Gamma_L \cap \Gamma_R}$ and the $QR$ factorization of $Y$ per iteration.\\

\begin{algorithm}\label{alg:7-20-1} (A FEAST algorithm for solving (\ref{equ:7-14-1}) with $\lambda \in (a, b)$)

Input:  $A \in {\mathbb R}^{n \times n}$ is symmetric, $Y \in {\mathbb R}^{n \times m}$ random with
$m \geq s$,
 $\Gamma_L$ and $\Gamma_R$ the circles described in \S\ref{sec:intro}, $max\_it$
a maximum number of iteration, $tol$ a convergence  tolerance.

Output: computed eigenvalues and eigenvectors are stored in $Eigvlu$ and $Eigvtr$ respectively.
\end{algorithm}

\begin{tabbing}
x\=xxx\= xxx\=xxx\=xxx\=xxx\=xxx\=xxx\=xxx\kill
\> Function $[ Eigvlu,  Eigvtr]$ = {\sc FEAST$_2$}$( A,  Y, \Gamma_L, \Gamma_R, max\_it, tol)$\\
\>1.\> For $iter = 1, \ldots, max\_it$\\
\>2.\>\> Compute $\displaystyle{Y = P_{\Gamma_L} Y}$ by (\ref{equ:7-15-2}).\\
\>3.\>\> Compute $\displaystyle{Y = P_{\Gamma_R} Y}$ by (\ref{equ:7-15-2}).\\
\>4.\>\> Compute the $QR$ factorization $Y = QR$
where $Q \in {\mathbb R}^{n \times m}$ and $R \in {\mathbb R}^{m \times m}$.\\
\>\>\> Set $Y = Q$.\\
\>5.\>\> Set $\hat{A}=  Y^T A  Y$ and solve the eigenproblem $\hat{A} \hat{x} = \hat{\lambda}\hat{x}$
to obtain the\\
\>\>\> eigenpairs $\{(\hat{\lambda}_i, \hat{x}_i)\}_{i = 1}^{m}$.\\
\>6.\>\> Compute $x_i = Y \hat{x}_i$ for $i =1,2,\ldots m$. \\
\>7.\>\> Calculate the maximum relative residual norm
$\tau = \max\{\|A x_i - \hat{\lambda}_i x_i\|_2/\|x_i\|_2,$\\
\>\>\> $1 \leq i \leq m, \hat{\lambda}_i \in (a, b)\}$. If $\tau < tol$, store the
eigenvalues $\hat{\lambda}_i \in (a, b)$ in $Eigvlu$\\
\>\>\> and their corresponding eigenvectors $x_i$ in $Eigvtr$, then break the for loop.\\
\>8.\> End
\end{tabbing}
\vspace{.2cm}

Theoretically the choice of the common radius $r$ of $\Gamma_L$ and $\Gamma_R$ is independent of the interval $(a, b)$ provided that $r \geq (b-a)/2$.
Computationally, however, $r$ should not be chosen arbitrarily large,
otherwise the subspace $span\{Y\}$ resulting from the computed $Y$ will be far from ${\mathbb K}$ due to computer rounding errors
and the truncation error of the quadrature rule in (\ref{equ:7-15-2}) and, as a result, the computed eigenpairs in Line 5 will not be accurate.

Algorithm \ref{alg:7-20-1}
requires $m \geq s$.
The restriction,
however, can be lifted by incorporating a power subspace iteration process into the algorithm.

\subsection{A power subspace iteration algorithm}\label{subsec: subiter} Power Subspace Iteration (PSI)
is an eigenvalue algorithm that permits us to compute a $m$-dimensional invariant subspace. It is
 a straightforward generalization of the power method for one eigenvector. Let us focus on real and symmetric matrices.
Given a symmetric $A \in {\mathbb R}^{n \times n}$
 and starting with $Y \in {\mathbb R}^{n \times m}$, the algorithm produces scalar sequences that approach the $m$ dominant eigenvalues of
 $A$ and vector sequences that approach the corresponding eigenvectors. The following is an implementation version of Algorithm \ref{alg:9-4-1}
 with $\rho(A) = A$ and the $QR$ factorization of $Y$ per iteration.\\

 \begin{algorithm}\label{alg:7-22-1} (A PSI algorithm)
 The input and output quantities $A$, $Y$, $max\_it$, $tol$, $Eigvlu$, and $Eigvtr$ are described in Algorithm \ref{alg:7-20-1} except that
$m$ does not need to be greater or equal to $s$.
\end{algorithm}

\begin{tabbing}
x\=xxx\= xxx\=xxx\=xxx\=xxx\=xxx\=xxx\=xxx\kill
\> Function $[Eigvlu,  Eigvtr]$ = {\sc PSI}$( A,  Y, max\_it, tol)$\\
\>1.\> For $iter = 1, \ldots, max\_it$\\
\>2.\>\> $QR$-factorize $Y = QR$ where $Q \in {\mathbb R}^{n \times m}$ and $R \in {\mathbb R}^{m \times m}$. Set $Y = Q$.\\
\>3.\>\> Set $\hat{A}=  Y^T A  Y$ and solve the eigenproblem $\hat{A} \hat{x} = \hat{\lambda}\hat{x}$
to obtain the\\
\>\>\> eigenpairs $\{(\hat{\lambda}_i, \hat{x}_i)\}_{i = 1}^{m}$.\\
\>4.\>\> Compute $x_i = Y \hat{x}_i$ for $i =1,2,\ldots m$. \\
\>5.\>\> Calculate the maximum relative residual norm
$\tau = \max\{\|A x_i - \hat{\lambda}_i x_i\|_2/\|x_i\|_2,$\\
\>\>\> $1 \leq i \leq m\}$. If $\tau < tol$, store the
eigenvalues $\hat{\lambda}_i$ in $Eigvlu$ and their corresponding\\
\>\>\> eigenvectors $x_i$ in $Eigvtr$, then break the for loop.\\
\>6.\>\> Set $Y = A Y$. \\
\>7.\> End
\end{tabbing}
\vspace{.2cm}

Assume that the eigenvalues of $A$ are arranged in deceasing order in size. That is,
$$
|\lambda_1|  \geq \ldots \geq |\lambda_m| > |\lambda_{m+1}| \geq \ldots |\lambda_n|.
$$
Then the rate of convergence of the $i$th computed eigenvector $x_i$ (i.e., the eigenvector
associated with $\lambda_i$) depends on the ratio $|\lambda_{m+1}/\lambda_i|$.
Precisely, the distance between $x_i$ at iteration $k$ and the true eigenvector $v_i$ associated with $\lambda_i$ is  $O(|\lambda_{m+1}/\lambda_i|^{k})$ (see, for instance, \cite{saad}).

\subsection{A FEAST$_2$-PSI algorithm}\label{subsec: FFS}
Power subspace iteration is usually used to find some eigenvalues of the largest magnitude in the spectrum
of a matrix $A$, but it
can also be used to find some eigenvalues of the largest magnitude in a given interval $(a, b)$.
In fact, when Algorithm \ref{alg:7-22-1} is applied to the matrix $A P_{\Gamma_K \cap \Gamma_R}$, or equivalently,
applied to $A$ but starting with $P_{\Gamma_L \cap \Gamma_R} Y$, the algorithm will converge to the first $m$
dominant eigenvalues in absolute value
in the interval $(a, b)$.
In this case, the $Y$ in
Algorithm \ref{alg:7-22-1} satisfies (\ref{equ:10-16-1}) in every iteration theoretically.

Satisfying the condition (\ref{equ:10-16-1}) is crucial in order to find eigenvalues in $(a, b)$. Computationally, however,
the columns of
$Y$ cannot strictly lie in ${\mathbb K}$ due to roundoff errors and the truncation error of a quadrature rule. As a result, Algorithm \ref{alg:7-22-1} will eventually converge to the dominant eigenvalues of the whole spectrum of $A$
rather than to the dominant eigenvalues in $(a, b)$. To avoid this happening, it is necessary to make a correction on
$Y$ by applying the operator $P_{\Gamma_L \cap \Gamma_R}$ on it from time to time during the iteration process of Algorithm \ref{alg:7-22-1}
in order that
(\ref{equ:10-16-1}) is kept satisfied as much as possible.

To speed up the convergence of Algorithm \ref{alg:7-22-1},
one can apply the algorithm to a shifted matrix $A - \sigma I$
with the shift number $\sigma$ being carefully chosen. The following Algorithm \ref{alg:7-28-1} is a refinement of
Algorithm \ref{alg:7-22-1}, in which we find the $m$ largest eigenvalues of $A$ in $(a, b)$ rather than find the $m$ eigenvalues with the largest magnitude
in $(a, b)$.

Let the eigenvalues of $A$ in $(a, b)$ be arranged decreasingly:
\begin{equation}\label{equ:11-2-1}
\lambda_1 \geq \ldots \geq \lambda_m > \lambda_{m+1} \geq \ldots \geq \lambda_s,
\end{equation}
and we want to find $\lambda_1, \ldots, \lambda_m$.
If Algorithm \ref{alg:7-22-1} is applied to  $A - \sigma I$,
the corresponding rate
of convergence for each computed eigenvector $x_i$
will depend on $|(\lambda_{m+1} - \sigma)/(\lambda_i - \sigma)|$.
Ideally, $\sigma$ is chosen to
minimize $\displaystyle{\max_{1 \leq i \leq m}|(\lambda_{m+1} - \sigma)/(\lambda_i - \sigma)|}$ while
($\lambda_i - \sigma$)'s are the dominant eigenvalues
in absolute value, so that
all the $m$ desired eigenvalues converge as fast as possible.
It is easy to see that the best possible
such a $\sigma$ is
\begin{equation}\label{equ:7-29-1}
\sigma = (\lambda_{m+1} + \lambda_s)/2 \approx (\lambda_{m+1} + a)/2.
\end{equation}

Note that the matrices $A$ and $A - \sigma I$ share the same eigenvectors.
When we apply Algorithm \ref{alg:7-22-1} to $A - \sigma I$ to obtain the $m$ largest eigenvalues $\lambda_1, \ldots, \lambda_m$
of $A$, the $A$ in Line 3 of the algorithm can be kept unchanged while the $A$ in Line 6 is replaced by $A - \sigma I$.

The previous algorithms use the same stopping criterion which relies on the relative residual $\|A x_i - \hat{\lambda}_i x_i\|_2/\|x_i\|_2$ of the computed eigenpair $(\hat{\lambda}_i, x_i)$. This stopping criterion is not good enough from our numerical experiments since the matrix $A$ is not scaled into a matrix of order one in magnitude. Instead, we shall adopt the following relative residual to monitor the convergence of $(\hat{\lambda}_i, x_i)$,
\begin{equation}\label{equ:9-26-18}
\frac{\|\frac{1}{\mu} A x_i - \frac{1}{\mu} \hat{\lambda}_i x_i\|_2}{\|x_i\|_2}
= \frac{\|A x_i - \hat{\lambda}_i x_i\|_2}{\mu \|x_i\|_2}
\end{equation}
where $\mu$ is a scale factor defined by
\begin{equation}\label{equ:7-30-1}
\mu = \sqrt{\frac{(Ay)^T (Ay)}{n}}
\end{equation}
where $y \in {\mathbb R}^n$ is a random vector with iid
elements from $N(0, 1)$,
the normal distribution with mean $0$ and variance $1$. Approximately $\mu$ is the square root of the average of the squares of the eigenvalues of $A$.

We now summarize the above discussions in the following Algorithms \ref{alg:7-28-1} and \ref{alg:7-29-1}.
Algorithm \ref{alg:7-28-1} is a revised version of Algorithm \ref{alg:7-22-1} applied to the matrix $A - \sigma I$. We have added in the algorithm some
tests to determine whether (\ref{equ:10-16-1}) is violated.
Instead of computing the $m$ dominant eigenvalues in absolute value
like Algorithm \ref{alg:7-22-1}, Algorithm \ref{alg:7-28-1} computes only
the first $num\_cmp$ largest eigenvalues of $A$ in $(a, b)$ where $num\_cmp$ is a positive integer not greater than $m$.\\

\begin{algorithm}\label{alg:7-28-1}
(A PSI algorithm for some largest eigenvalues in $(a, b)$)

Input: $A \in {\mathbb R}^{n \times n}$ is symmetric and $Y \in {\mathbb R}^{n \times m}$ satisfies
(\ref{equ:10-16-1}). $num\_cmp$ is a positive integer not greater than $m$,
$num\_eigm$ a positive integer,
$\mu$ the scale factor in (\ref{equ:7-30-1}), $\theta$ a real number in the interval $[0, m]$, $max\_it$ a maximum number of iteration.
$EigmHist$ holds ``$num\_eigm$'' of the most recent estimates
of the $\lambda_m$ in (\ref{equ:11-2-1}).

Output: computed eigenvalues are stored in $Eigvlu$ in decreasing order, their corresponding eigenvectors in $Eigvtr$, and
their corresponding relative residual norms defined in (\ref{equ:9-26-18}) in $ErrList$. $Y \in {\cal R}^{n \times m}$ is the iteration matrix,
$iter$ the number of iterations performed.
$EigmHist$ holds ``$num\_eigm$'' of
the most recent estimates of the $\lambda_m$ in (\ref{equ:11-2-1}).
\end{algorithm}


\begin{tabbing}
x\=xxx\= xxx\=xxx\=xxx\=xxx\=xxx\=xxx\=xxx\=xxx\=xxx\kill
\>Function
$[Eigvlu, Eigvtr, ErrList, Y, iter, EigmHist]$\\
\> $=$ {\sc
PSI}$(A, Y, num\_cmp, num\_eigm, a, b,
\mu, \theta, EigmHist, max\_it)$\\
\>1.\> Set $count0 = -1, err0 = \infty$.\\
\>2.\>For $iter = 1, \dots, max\_it$\\
\>3.\>\>$QR$-factorize $Y = QR$ where $Q \in {\mathbb R}^{n \times m}$ and $R \in {\mathbb R}^{m \times m}$. Set $Y = Q$.\\
\>4.\>\> Set $\hat{A} = Y^TAY$ and solve the eigenproblem $\hat{A} \hat{x} = \hat{\lambda} \hat{x}$ to obtain the \\
\>\>\> eigenpairs $\{(\hat{\lambda}_i, \hat{x}_i)\}_{i = 1}^m$.\\
\>5.\>\>Compute $x_i = Y \hat{x}_i$ for $i = 1, 2, \ldots, m$.\\
\>6.\>\> Set $Eigvlu = [\,\,]$, $Eigvtr = [\,\,]$, and $ErrList = [\,\,]$.\\
\>7.\>\> Set $err = -1$, $count = 0$, and $count1 = 0$. \,\,\,\,\,\,\, \% $count1$ is an estimate of $\min\{m, s\}$.\\
\>8.\>\> Set $eigm = \infty$. \,\,\,\,\,\,\, \% $eigm$ is an estimate of $\lambda_m$ in (\ref{equ:11-2-1}).\\
\>9.\>\> For $i = 1, \ldots, m$\\
\>10.\>\>\> Determine $i_0$ so that $\hat{\lambda}_{i_0} = \max\{\hat{\lambda}_1, \hat{\lambda}_2, \ldots,
\hat{\lambda}_m\}$.\\
\>11.\>\>\> If $\hat{\lambda}_{i_0} \in (a, b)$\\
\>12.\>\>\>\> $count1 = count1 + 1$.\\
\>13.\>\>\>\> If $count1 \leq num\_cmp$\\
\>14.\>\>\>\>\> Compute $erri = \|Ax_{i_0} - \hat{\lambda}_{i_0}x_{i_0}\|_2/(\mu \|x_{i_0}\|_2)$.\\
\>15.\>\>\>\>\> $count = count + 1$;\\
\>16.\>\>\>\>\> $Eigvlu = [Eigvlu, \hat{\lambda}_{i_0}]$; $Eigvtr = [Eigvtr, x_{i_0}]$;\\
\>17.\>\>\>\>\> $ErrList = [ErrList, erri]$;\\
\>18.\>\>\>\>\> $err = \texttt{max}(err, erri)$.\\
\>19.\>\>\>\>End\\
\>20.\>\>\>\>$eigm = \texttt{min}(eigm, \hat{\lambda}_{i_0})$.\\
\>21.\>\>\>End\\
\>22.\>\>\>Set $\hat{\lambda}_{i_0} = a - 1$. \\
\>23.\>\>End\\
\>24.\>\>\\
\>25.\>\>If ($count = count0$) and $(err > err0)$\\
\>26.\>\>\>$Eigvlu = Eigvlu0$; $Eigvtr = Eigvtr0$;\\
\>27.\>\>\>$ErrList = ErrList0$; $EigmHist = EigmHist0$;\\
\>28.\>\>\>$Y = Y0$;
$iter = iter - 1$.\\
\>29.\>\>\>Break the $iter$-for loop.\\
\>30.\>\>End\\
\>31.\>\>\\
\>32.\>\>If $(count \ne count0) \,\, \& \,\, (iter > 1)$\\
\>33.\>\>\>$Eigvlu = Eigvlu0$; $Eigvtr = Eigvtr0$;\\
\>34.\>\>\>$ErrList = ErrList0$; $EigmHist = EigmHist0$;\\
\>35.\>\>\>$Y = Y0$;
$iter = iter - 1$.\\
\>36.\>\>\>Break the $iter$-for loop.\\
\>37.\>\>End\\
\>38.\>\>\\
\>39.\>\> If $(m > \theta \cdot count1)$ or $(max\_it = 1)$\\
\>40.\>\>\>Break the $iter$-for loop.\\
\>41.\>\>End\\
\>42.\>\>\\
\>43.\>\>$Eigvlu0 = Eigvlu$; $Eigvtr0 = Eigvtr$; \,\,\, \% keep the data in current iteration.\\
\>44.\>\>$ErrList0 = ErrList$; $EigmHist0 = EigmHist$;\\
\>45.\>\>$Y0 = Y$; $err0 = err$; $count0 = count$.\\
\>46.\>\>\\
\>47.\>\>$EigmHist = [EigmHist, eigm]$. \\
\>48.\>\>$leng = \texttt{length}(EigmHist)$.\\
\>49.\>\>If $leng > num\_eigm$\\
\>50.\>\>\>$EigmHist = EigmHist(2:leng)$; \,\,\,\,\,\,\,\% $EigmHist$ keeps ``$num\_eigm$'' of \\
\>\>\>\>\>\>\>\>\>\> \% the most recent estimates of $\lambda_m$ in (\ref{equ:11-2-1}).\\
\>51.\>\>\>$leng = leng - 1$.\\
\>52.\>\>End\\
\>53.\>\>$eigm = \texttt{sum}(EigmHist)/leng$. \, \% taking an average
makes $eigm$ less vibrating.\\
\>54.\>\>$\sigma = (eigm + a)/2$;\\
\>55.\>\>$Y = (A - \sigma I)Y$.\\
\>56.\>End
\end{tabbing}
\vspace{.2cm}

In Lines 6-23, Algorithm \ref{alg:7-28-1} selects the first ``$num\_cmp$'' largest
eigenvalues $\hat{\lambda}_i$ that lie in the interval $(a, b)$ and their corresponding
eigenvectors $x_i$.
In Lines 25-41,
several tests for the violation of condition (\ref{equ:10-16-1}) are presented. The design of the violation tests
is
similar to that of
the stopping criteria in Algorithm 5 of \cite{ycy}, and
works well in our numerical experiments.
When a violation test is passed, the $iter$-for loop is stopped, and the algorithm outputs the
iteration matrix $Y$ for a correction. The following Algorithm \ref{alg:7-29-1}
will perform the correction (in its Lines 5 and 6) by pre-multiplying $Y$ with $P_{\Gamma_L \cap \Gamma_R}$.

About Line 39, we consider $count1$ as an estimate of $\min\{m, s\}$. So ``$m > \theta \cdot count1$'' can be understood as ``$m > \theta \min\{m, s\}$''.
When ``$m > \theta \min\{m, s\}$'' is true, we break the $iter$-for loop and do not perform the following Lines 43 - 56 which are about the power subspace iteration.

There are two extreme values for $\theta$: (i) $\theta = 0$, Lines 43 - 56 are never implemented; (ii) $\theta = m$, Lines 43 - 56 are always implemented unless $count1 = 0$ or $max\_it = 1$.

Ideally the shift $\sigma$ is computed by (\ref{equ:7-29-1}), but this is impossible because the information about $\lambda_{m+1}$ is not available in the algorithm.
Instead, we use the equation
\begin{equation}\label{equ:9-30-1}
\sigma = (\lambda_m + a)/2
\end{equation}
to compute $\sigma$ in Line 54 with $\lambda_m$ estimated in Line 53.

Our experiments showed that the eigenvalues computed by Algorithm \ref{alg:7-28-1} converge at different rates, usually with the larger ones converging faster.
So, instead of outputting all the ``$num\_cmp$'' computed eigenvalues, Algorithm \ref{alg:7-29-1} below only outputs those with higher convergence
rates. More precisely, Algorithm \ref{alg:7-29-1} outputs ``$num\_out$'' of the ``$num\_cmp$'' eigenvalues
computed by Algorithm \ref{alg:7-28-1} where ``$num\_out$'' is a positive integer not greater than ``$num\_cmp$''.
\\

\begin{algorithm}\label{alg:7-29-1}
(A FEAST$_2$-PSI algorithm for some largest eigenvalues in $(a, b)$)

Input: $A \in {\mathbb R}^{n \times n}$ is symmetric, $Y \in {\cal R}^{n \times m}$ random, $\Gamma_L$ and $\Gamma_R$
the circles in \S\ref{sec:intro}.
The quantities $\theta$,
$num\_cmp$, $num\_eigm$,
and $max\_it$ are described in Algorithm \ref{alg:7-28-1}.
$num\_out$, a positive integer not greater than ``$num\_cmp$'', is the number
of output eigenpairs, and
$sub\_max\_it$ a maximum number of iteration used by Algorithm \ref{alg:7-28-1}.

Output: the first ``$num\_out$'' largest eigenvalues of $A$ in $(a, b)$ are
output and stored in decreasing order in $Eigvlu$,
their corresponding eigenvectors in $Eigvtr$,
and their corresponding relative residual norms in $ErrList$. The output eigenpairs have the smallest
maximum-relative-residual-norm.
$ErrHist$ keeps the history of maximum-relative-residual-norm per iteration, and $NumAYHist$ the history of the number of
$A - \sigma I$ times $Y$
per iteration.
\end{algorithm}

\begin{tabbing}
x\=xxx\= xxx\=xxx\=xxx\=xxx\=xxx\=xxx\=xxx\=xxx\=xxx\=xxx\=xxx\=xxx\kill
\>Function $[Eigvlu, Eigvtr, ErrList, ErrHist, NumAYHist]$\\
\>$=$ {\sc F$_2$P}$(A, Y, \Gamma_L, \Gamma_R, a, b, \theta, num\_cmp, num\_out, num\_eigm,
max\_it, sub\_max\_it)$\\
\>1.\>Set $Eigvlu = [\,\,]$, $Eigvtr = [\,\,]$, and $ErrList = [\,\,]$. \\
\>2.\>Set $ErrHist = [\,\,]$, $NumAYHist = [\,\,]$, and $EigmHist = [\,\,]$.\\
\>3.\> $count = 0$, $err0 = \infty$ and compute the scale factor $\mu$ according to (\ref{equ:7-30-1}).\\
\>4.\>For $iter = 1, \ldots, max\_it$\\
\>5.\>\> Compute $\displaystyle{Y = P_{\Gamma_L} Y
}$ by (\ref{equ:7-15-2}).\\
\>6.\>\> Compute $\displaystyle{Y = P_{\Gamma_R} Y
}$ by (\ref{equ:7-15-2}).\\
\>7.\>\>$[Eigvlu1, Eigvtr1, ErrList1, Y, sub\_iter, EigmHist]$\\
\>\>\> $=$ {\sc PSI
}$(A, Y, num\_cmp, num\_eigm, a, b,
\mu, \theta, EigmHist, sub\_max\_it)$;\\
\>8.\>\>$NumAYHist = [NumAYHist, sub\_iter - 1]$; \,\, \% $sub\_iter - 1$ is the number\\
\>\>\>\>\>\>\>\>\>\>\>\>\,\,\,\,\,\,\,\,\,\,\,\,\,\,\,\,\, \% of $A - \sigma I$ times $Y$ performed by PSI\\
\>9.\>\\
\>10.\>\> $err = -1$; $leng = \texttt{min}(num\_out, \texttt{length}(Eigvlu1))$.\\
\>11.\>\>If $leng > 0$\\
\>12.\>\>\>$count = count + 1$.\\
\>13.\>\>\>$err = \texttt{max}(ErrList1(1: leng))$. \,\,\,\% $err$ is a maximum-relative-residual-norm\\
\>14.\>\>End\\
\>15.\>\>$ErrHist = [ErrHist, err]$.\\
\>16.\>\\
\>17.\>\>If $err < err0$ and $err \ne -1$\\
\>18.\>\>\>$Eigvlu = Eigvlu1(1: leng)$; $Eigvtr = Eigvtr1(:, 1: leng)$;\\
\>19.\>\>\>$ErrList = ErrList1(1: leng)$; $err0 = err$.\\
\>20.\>\>End\\
\>21.\>End\\
\>22.\>If $count < max\_it - count$ \,\,\, \%\,\,$count$ is the number of $err \ne -1$.\\
\>23.\>\> $Eigvlu = [\,\,], Eigvtr = [\,\,]$, and $ErrList = [\,\,]$.\\
\>24.\>End\\
\end{tabbing}


Algorithm \ref{alg:7-29-1} is a combination of a spectrum projection subspace iteration (Lines 5-6) and a power subspace iteration (Line 7).
It uses the spectrum projection iteration to keep the computed $Y$ satisfying (\ref{equ:10-16-1}) and the power iteration to find the desired eigenvalues and eigenvectors.

We remark that Algorithm \ref{alg:7-29-1} will reduce to a version of Algorithm \ref{alg:7-20-1} with a different stopping criterion if we set $m \geq s$
and $sub\_max\_it = 1$ in the algorithm. In the case when $sub\_max\_it = 1$, Algorithm \ref{alg:7-29-1} does not involve any power subspace iteration.

\subsection{All the eigenvalues in an interval}\label{subsec:alleigs}
Unlike the FEAST algorithm, Algorithm \ref{alg:7-29-1}
only finds some of the largest eigenvalues of $A$ in a given interval $(a, b)$. If one wants to find all the eigenvalues in $(a, b)$ by Algorithm \ref{alg:7-29-1},
here is a strategy for achieving the goal. First, apply Algorithm \ref{alg:7-29-1} to the interval $(a, b)$ to obtain some largest eigenvalues $\hat{\lambda}^{(1)}_1 \geq \hat{\lambda}^{(1)}_2 \geq \ldots \geq \hat{\lambda}^{(1)}_{k_1}$ in $(a, b)$. When this is done, pick a point $b_1$ between $\hat{\lambda}^{(1)}_1$ and $\hat{\lambda}^{(1)}_{k_1}$ and set $a_1 = b_1 - \delta$, where $\delta = b - a$.
Then, apply the algorithm to the interval $(a_1, b_1)$ to obtain some largest eigenvalues $\hat{\lambda}^{(2)}_1 \geq \hat{\lambda}^{(2)}_2 \geq \ldots \geq \hat{\lambda}^{(2)}_{k_2}$ in $(a_1, b_1)$.
When this is done, pick a point $b_2$ between $\hat{\lambda}^{(2)}_1$ and $\hat{\lambda}^{(2)}_{k_2}$ and set $a_2 = b_2 - \delta$,
then apply the algorithm to the interval $(a_2, b_2)$. This process is repeated until all the eigenvalues in $(a, b)$ have been found.

If a given interval
is large, one can divide it into smaller subintervals, then apply the strategy in parallel to each of the subintervals.

\section{Numerical Experiments}\label{sec:numexp} In this section, we present some
experiments to illustrate the behavior
of Algorithm \ref{alg:7-29-1} using three test matrices.
Two of them are from The University of Florida Sparse
Matrix Collection \cite{ufl_sparse_matrics} described below:
\begin{enumerate}
\item[(1)]
{\it Na5} is a $5832 \times 5832$ real and symmetric matrix with $305,630$ nonzero entries, from a theoretical/quantum chemistry problem. The spectrum range of the matrix is $[-0.1638, 25.67]$
and the Average Number of Eigenvalues (ANE) is $5832/(25.67 - (-0.1638)) \approx 225.75$ eigenvalues per unit interval.
\item[(2)]
{\it Andrews} is a $60000 \times 60000$ real and symmetric matrix with $760,154$ nonzero integer entries, from a computer graphics/vision problem. The spectrum range is $[0, 36.49]$ and
 the ANE is $60000/(36.49 - 0) \approx 1644.29$ eigenvalues per unit interval.
\end{enumerate}
The third test matrix is a random diagonal matrix of size $10^6$.

All the computations are carried out in \textsc{Matlab} version R2017b on a Windows 10 machine.
The eigenproblem $\hat{A}\hat{x} = \hat{\lambda} \hat{x}$ in Line 4 of Algorithm \ref{alg:7-28-1} is solved by the \textsc{Matlab} function \texttt{eig}.
Except for the third test matrix, in numerical comparisons, we treat the eigenvalues and eigenvectors computed by \texttt{eig} or \texttt{eigs} as the exact eigenvalues and eigenvectors,
and results obtained by Algorithm \ref{alg:7-29-1} are compared to them\footnote{The eigenvalues and eigenvectors computed by \texttt{eig} or \texttt{eigs} may not necessarily be accurate. Our experiments show that the eigenpairs computed by Algorithm \ref{alg:7-29-1} are at the same level of accuracy with those obtained by \texttt{eig} or \texttt{eigs}.}.
The construction of the third test matrix ensures that the eigenvalues are known to us and are the diagonal entries of the matrix.

We use the Gauss-Legendre quadrature rule on the interval $[0, 1]$ with $q = 8$ in (\ref{equ:7-15-2}).
As for the solution of the $mq$ linear systems in (\ref{equ:7-15-3}), we employ the two-term recurrence Krylov subspace
solver BiCG\cite{Fletcher}. BiCG requires two matrix-vector multiplications per iteration and is robust in performance.
We solve
the linear systems 
sequentially with initial guess $x = 0$ and the stopping criterion
$\|r\|_2/\|b\| < 10^{-10}$, where $b$ represents the right hand side of a linear system and $r$ the computed residual vector.
We remark that BiCG can be replaced by any other linear solver (see, for instance, \cite{gvl,saad03, Vorst} for other linear solvers), and
we also note that when $A$ is real and symmetric,
the MINRES method\cite{ps}, a symmetric version of GMRES\cite{saad03, gmres}, is a good choice since
it requires only one matrix-vector multiplication per iteration.

The following values of the input arguments of Algorithm \ref{alg:7-29-1} are fixed
for all the experiments:
$Y \in {\mathbb R}^{n \times m}$ is a random matrix with iid elements from $N(0, 1)$, $\theta = 3$, $num\_eigm = 5$, $max\_it = 50$ or $100$,
$num\_cmp = \lfloor m/2 \rfloor$,
$num\_out = \lfloor m/2 \rfloor$ or $\lfloor m/4 \rfloor$
where $\lfloor \,
\cdot \,\rfloor$ rounds its argument to the nearest integer towards
minus infinity. The common radius $r$ of the circles $\Gamma_L$ and $\Gamma_R$ is set to $5$ for
{\it Na5}, and to $2$ for
{\it Andrews} and the random diagonal matrix in \S\ref{sec:10-2-1}.
Further,
$sub\_max\_it = 100$ except otherwise specified.


\subsection{Experiments with Na5} The eigenvalues $\lambda_i$ and eigenvectors $v_i$ of the matrix {\it Na5} computed by
\texttt{eig} satisfy $\displaystyle{
\max_{1 \leq i \leq n} \|Av_i - \lambda_i v_i\|_2/\|v_i\|_2}$
$< 1.02 \times 10^{-13}
$ where $n = 5832$ is the size of {\it Na5}. For this matrix, the scale factor $\mu$ in (\ref{equ:7-30-1}) is
about $11.72$.\\

{\bf Experiment 1.}
In this experiment, we compare the performances of FEAST, FEAST$_2$, and F$_2$P in terms of the minimum maximum
relative residual norm $\tau_r$ and the maximum relative error $\tau_\lambda$ defined as follows.
\begin{equation}\label{equ:8-10-1}
\tau_r = \min_k \tau_r^{(k)} \equiv \min_{k} \max_i \| A x_i^{(k)} - \hat{\lambda}_i^{(k)} x_i^{(k)} \|_2/(\mu \|x_i^{(k)}\|_2)
\end{equation}
where $(\hat{\lambda}_i^{(k)}, x_i^{(k)})$'s are the eigenpairs computed by Algorithm \ref{alg:7-29-1} in its $k$th iteration.\footnote{The eigenpairs $(\hat{\lambda}_i^{(k)}, x_i^{(k)})$ are those output by the function PSI in Line 7 of Algorithm \ref{alg:7-29-1}.}
Moreover, let $k_0$ be the iteration number that satisfies $\displaystyle{\tau_r^{(k_0)} = \min_k \tau^{(k)}_r}$. We then define
\begin{equation}\label{equ:8-10-2}
\tau_\lambda = \max_i | \hat{\lambda}_i^{(k_0)} - \lambda_i |/|\lambda_i|.
\end{equation}
If one makes use of the output quantities of the function $F_2P$ in Algorithm \ref{alg:7-29-1},
then $\tau_r = \texttt{max}(ErrList)$, the entries of the vector $ErrHist$ are the $\tau_r^{(k)}$'s, and the entries of
$Eigvlu$ are the $\hat{\lambda}_i^{(k_0)}$'s.


In this experiment, instead of using Algorithm \ref{alg:9-4-1} with $\rho(A) = Q_{(a, b)} Q_{(a, b)}^T$ as our FEAST algorithm, we use the one described here: in Algorithm \ref{alg:7-29-1}, we delete Line 6, and set $\Gamma_L$ to be the circle with center $c = (a+b)/2$ and radius $r = (b - a)/2$;
in addition, we set $sub\_max\_it = 1$. Similarly, we employ Algorithm \ref{alg:7-29-1} with
$sub\_max\_it = 1$ rather than Algorithm \ref{alg:7-20-1} as our FEAST$_2$ algorithm. For F$_2$P algorithm, we use Algorithm \ref{alg:7-29-1}.

\begin{figure}
\begin{center}
\includegraphics[width=6.4cm]{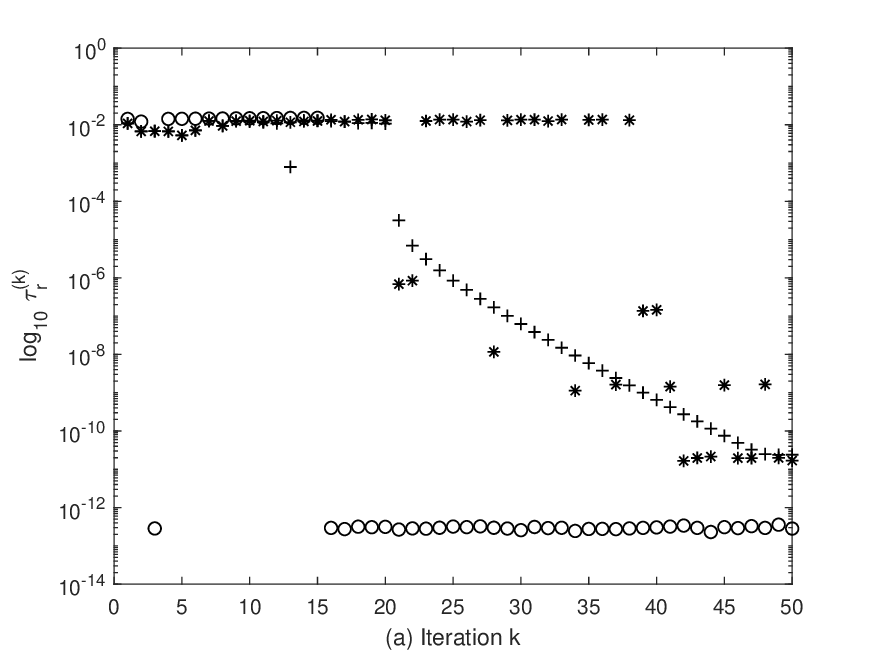}
\includegraphics[width=6.4cm]{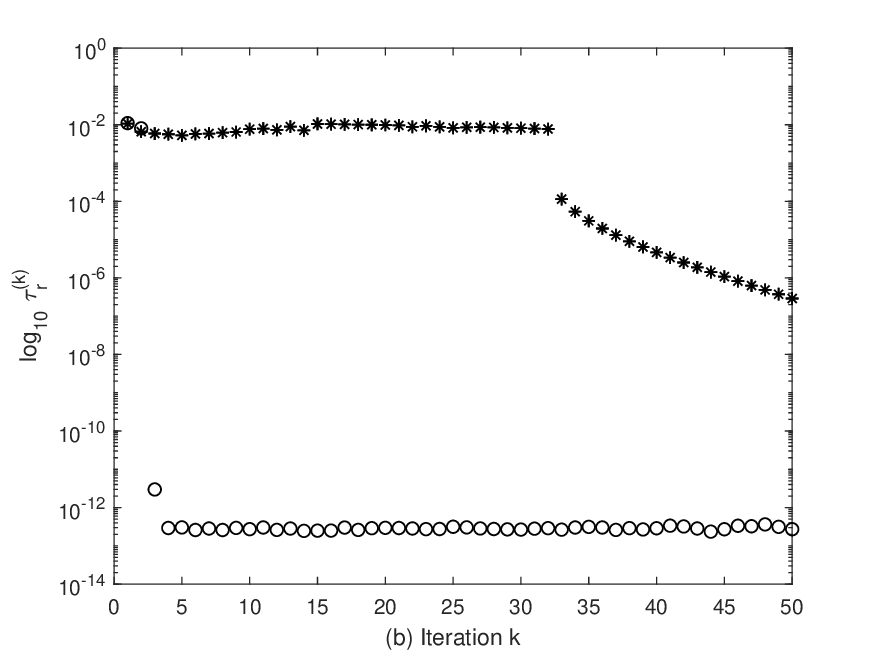}
\caption{Experiment 1: convergence histories of the FEAST, FEAST$_2$, and F$_2$P algorithms. FEAST: $\circ$; FEAST$_2$: $+$; F$_2$P: $*$. (a) $m = 130$, $num\_cmp = 65$, and $num\_out = 32$. (b) $m = 110$, $num\_cmp = 55$, and $num\_out = 27$.}
\label{Fig:8-15-1}
\end{center}
\end{figure}

\begin{figure}
\begin{center}
\includegraphics[width=6.4cm]{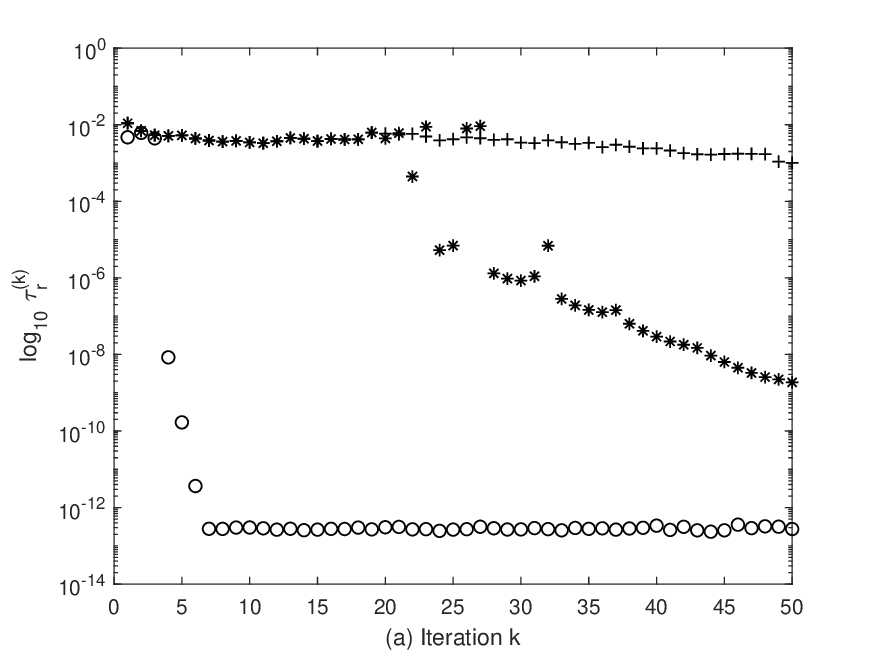}
\includegraphics[width=6.4cm]{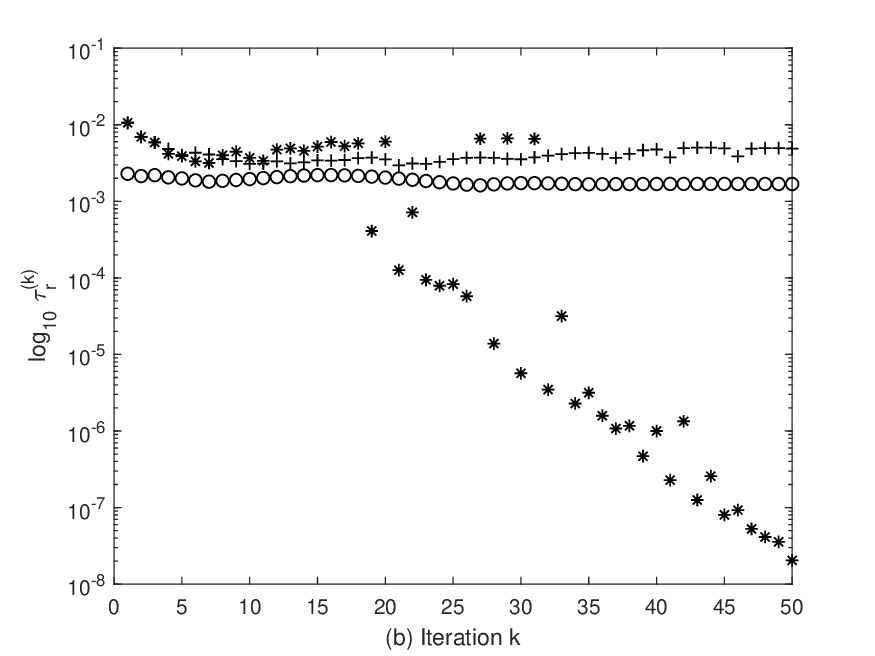}
\caption{Experiment 1: convergence histories of the FEAST, FEAST$_2$, and F$_2$P algorithms. FEAST: $\circ$; FEAST$_2$: $+$; F$_2$P: $*$. (a)
 $m = 90$, $num\_cmp = 45$, and $num\_out = 22$. (b) $m = 70$, $num\_cmp = 35$, and $num\_out = 17$.}
\label{Fig:9-4-1}
\end{center}
\end{figure}

Consider now the interval $(a, b) = (11.8, 12)$. This interval contains $84$ eigenvalues of $A$
and is in the middle of the spectrum.
Numerical results are presented in Tables \ref{Tab:8-6-1} and \ref{Tab:8-30-1}. The only different setting in the tables is that we set
$num\_out = \lfloor m/2 \rfloor$ in Table \ref{Tab:8-6-1} and $num\_out = \lfloor m/4 \rfloor$ in Table \ref{Tab:8-30-1}.
We also plot $\tau_r^{(k)}$ against the iteration number $k$ for each case in
Table \ref{Tab:8-30-1} in Figures \ref{Fig:8-15-1} and \ref{Fig:9-4-1}.

From Tables \ref{Tab:8-6-1} and \ref{Tab:8-30-1}, we observe that the computed eigenpairs have different convergence rates, with those associated with
more dominant eigenvalues converging faster.
In fact, the eigenpairs associated with the first $\lfloor m/4 \rfloor$ largest eigenvalues
converge faster than those associated with the first $\lfloor m/2 \rfloor$ largest eigenvalues since the corresponding $\tau_r$'s and $\tau_\lambda$'s
in Table \ref{Tab:8-30-1} are smaller.

FEAST$_2$ is a FEAST algorithm employing
two contour integrals per iteration.
By comparison,
FEAST has a more robust performance in terms of smaller $\tau_r$'s and $\tau_\lambda$'s. Note that
the circle radius chosen for
FEAST is $r = 0.1$ and that for FEAST$_2$ is $r = 5$. We believe the more robustness of FEAST is due to a smaller circle radius.
A side effect of a smaller circle radius, however, is that BiCG takes more iterations to converge in the solution of the linear systems in (\ref{equ:7-15-3}).
It is because those systems
may become more ill-conditioned when a smaller
circle radius is used (see (\ref{equ:7-18-1})).

In the case when $m = 70$, both FEAST and FEAST$_2$ fail to converge due to the violation of the
necessary condition $m \geq s$.
On the other hand, F$_2$P, an algorithm of FEAST$_2$ plus a power
iteration process, converges well. This experiment demonstrates the usefulness of adding a power iteration process to a FEAST algorithm.


In another case when $m = 110$, F$_2$P and FEAST$_2$
are identical since $\#(A - \sigma I)Y = 0$, namely, the power iteration process is not implemented.\\


\begin{table}
\centering
\caption{Experiment 1:
$m$ is the column size of the iteration matrix $Y$, $nc (= num\_cmp)$ and $no(= num\_out)$ the input arguments of Algorithm \ref{alg:7-29-1},
$\tau_r$ and $\tau_\lambda$ defined by (\ref{equ:8-10-1}) and (\ref{equ:8-10-2}),
$\#$iter the
maximum number of BiCG iterations required by the linear system in (\ref{equ:7-15-3}) that takes the longest to converge,
and $\#(A - \sigma I)Y$ the total number of the shifted matrix
$(A - \sigma I)$ times $Y$ performed by F$_2$P, i.e. $\#(A - \sigma I)Y = \texttt{sum}(NumAYHist)$ where $NumAYHist$ is the
output quantity of Algorithm \ref{alg:7-29-1}.
In this experiment, we set $max\_it = 50$. Let $\#$eig\_out $= \texttt{length}(Eigvlu)$ where $Eigvlu$ is the
output quantity of Algorithm \ref{alg:7-29-1}.
In each case of this experiment, we observed that
$\#$eig\_out was equal to $num\_out$.
}
\footnotesize{
\noindent
\begin{tabular}{|ccc|cccc|ccc|} \hline
 & & & &  FEAST&  &  &  &FEAST$_2$ &  \\ \hline
$m$ & nc & no &$\tau_r$ &$\tau_\lambda$ &$\#$iter&  &$\tau_r$ &$\tau_\lambda$&$\#$iter\\ \hline
130 &65 & 65 &1.39e-02&5.43e-04 &15761 &  & 6.79e-03&1.30e-03 &970 \\ \hline
110  &55 & 55 &2.34e-13 &5.18e-15 &15761  &  &2.90e-07 & 3.07e-12 & 962\\ \hline
90 & 45& 45 & 2.34e-13 &5.52e-15 &15761  &  &3.63e-03 &6.36e-04 &962 \\ \hline
70 & 35& 35 & 2.13e-03&1.52e-03 &15777  &  &3.08e-03 &2.06e-03 & 970\\ \hline\hline
 & & & &  F$_2$P& & & & & \\ \hline
$m$ & nc & no  & $\tau_r$& $\tau_\lambda$ &$\#$iter &$\#(A - \sigma I)Y$ &  & &  \\ \hline
130 &65 & 65&1.71e-11 &4.17e-15 & 970 &14 & & & \\ \hline
110 &55 &55 &2.90e-07 & 3.07e-12&962 &0 & & & \\ \hline
90 &45 &45 & 2.23e-09 &3.43e-15 &962 &10 & & & \\ \hline
70 & 35&35 & 5.18e-05 &4.17e-08 &970 & 28 & & & \\ \hline
\end{tabular}}
\label{Tab:8-6-1}
\end{table}

\begin{table}
\centering
\caption{
Experiment 1: for the meanings of the quantities $m$, nc, no, $\tau_r$, $\tau_\lambda$,
$\#$iter, $\#(A - \sigma I)Y$, and $\#$eig\_out, refer to the caption of Table \ref{Tab:8-6-1}.
In this experiment, we set $max\_it = 50$.
In each case of this experiment, we observed that $\#$eig\_out $= num\_out$.
}
\footnotesize{
\noindent
\begin{tabular}{|ccc|cccc|ccc|} \hline
 & & & & FEAST  & & & &FEAST$_2$ &  \\ \hline
$m$ & nc & no &$\tau_r$ &$\tau_\lambda$ &$\#$iter& &  $\tau_r$ &$\tau_\lambda$&$\#$iter \\ \hline
130 &65 & 32 &2.30e-13 &6.25e-15 &15761  & &2.37e-11  & 5.64e-15& 970 \\ \hline
 110 & 55 & 27& 2.35e-13&5.18e-15 &15761& &2.90e-07 &3.07e-12 &962  \\ \hline
 90 & 45 & 22 &2.34e-13 &4.75e-15 &15761 & &1.00e-03&2.76e-05&962\\ \hline
 70 & 35 & 17 & 1.62e-03&1.27e-03 &15777 & &2.94e-03 &2.13e-03 &970 \\ \hline\hline
 & & & &F$_2$P& & & & & \\ \hline
$m$ & nc & no & $\tau_r$& $\tau_\lambda$ &$\#$iter &$\#(A - \sigma I)Y$ & & & \\ \hline
130 & 65& 32&1.65e-11 &6.54e-15 &970 & 14& & &  \\ \hline
 110 & 55 & 27&2.90e-07 &3.07e-12& 962&0 & & & \\ \hline
 90 & 45 & 22& 1.86e-09&6.24e-15 &962 &10 & & &\\ \hline
 70 & 35 & 17 & 2.04e-08& 1.39e-14& 970 & 28 & & &\\ \hline
\end{tabular}}
\label{Tab:8-30-1}
\end{table}

{\bf Experiment 2.} In \S\ref{subsec: FFS}, we have seen that the convergence rate of the $i$th eigenvector $x_i$ computed by
Algorithm \ref{alg:7-28-1}
depends on the ratio $|(\lambda_{m+1} - \sigma)/(\lambda_i - \sigma)|$ where $\sigma$ is
given by (\ref{equ:9-30-1}). There are two situations in which this ratio is likely to be close to $1$, and as a result the
convergence may be slow: (i) $\lambda_i$ and $\lambda_{m+1}$ are likely to be close to each other when
the exact number $s$ of eigenvalues in the interval $(a, b)$ is much larger than $m$; (ii) the shift
$\sigma$ is likely to be far from both $\lambda_i$ and $\lambda_{m+1}$ when the interval $(a, b)$ is large. In this experiment, we demonstrate the
behaviors of Algorithm \ref{alg:7-29-1} in the two situations. We also show the ability of
Algorithm \ref{alg:7-29-1} catching eigenvalues when the interval $(a, b)$ is narrow relative to $\Delta$.
We remark that situations where the spectrum of $A$ may cause slow or varying convergence rates of FEAST are discussed in \cite{GP16, GPTV, xisaad}
and relevant remedies are provided there.

We now pick for $(a, b)$ a sequence of intervals in decreasing length, and we fix $m = 60$, $num\_cmp = 30$, $num\_out = 15$
in Algorithm \ref{alg:7-29-1}. About $max\_it$, we first set it to $50$, then increase it to $100$. The numerical results are listed in Tables \ref{Tab:9-11-1} and \ref{Tab:9-14-1}. From the tables,
we can see that Algorithm \ref{alg:7-29-1} converges
slowly when $(a, b) = (11.5, 12)$. The corresponding $s = 200$ which is much larger than $m$. As we reduce the length of the interval $(a, b)$, however, the number $s$ of eigenvalues in $(a, b)$ decreases accordingly and the algorithm tends to converge faster.
Moreover, it is interesting to see that the algorithm is
capable of catching the eigenvalues in $(a, b)$
accurately even when $(a, b)$ is very narrow, given that the equal radii $r$ of the circles $\Gamma_L$ and $\Gamma_R$ is $5$, a relatively large number.

We also observe that, even though it is not required to be greater than $s$ in Algorithm \ref{alg:7-29-1}, $m$ loosely depends on $s$ computationally. It should be
chosen near $s$ in order that the algorithm behaves well. Techniques of efficiently estimating the value of $s$ have been developed in the literature, see, for instance, \cite{futa, kpt, NPY, TP13, yxccb}.
Moreover, to reduce the dependence of $m$ on $s$,
a spectral transformation \cite{Parlett}, in particular, a transformation made by a certain
polynomial, may be needed \cite{bks, FangSaad}.

The following two phenomena about Algorithm \ref{alg:7-29-1} are also observed in this experiment. First, the computed eigenvalues seem to converge
faster than their associated computed eigenvectors since $\tau_\lambda$'s are generally smaller than their corresponding $\tau_r$'s by several orders of magnitude.
To speed up the convergence of the computed eigenvectors,
one idea may be adaptively
decreasing the common radius $r$ of the circles $\Gamma_L$ and $\Gamma_R$.
Moreover, using the filters
and techniques
introduced in \cite{GPTV, kbn, bib:WD2019, xisaad} may also be helpful.
Second, in the case when $s > m$,
the total number
$\#(A - \sigma I)Y$
of $(A - \sigma I)$ times $Y$ performed by the algorithm is considerably large. This hints that the power
subspace iteration plays a heavy role in the convergence of the algorithm.
When $s < num\_cmp$,
on the other hand,
the spectrum projection subspace iteration is dominant since $\#(A - \sigma I)Y \approx 0$.
We remark that
Algorithm \ref{alg:7-29-1} is reduced to a FEAST$_2$ algorithm in the case when $\#(A - \sigma I)Y = 0$.\\

\begin{table}
\centering
\caption{Experiment 2:
Results about Algorithm \ref{alg:7-29-1}. We set
$m = 60$, $num\_cmp = 30$, $num\_out = 15$, and $max\_it = 50$.
For the meanings of the quantities $\tau_r$, $\tau_\lambda$,
$\#$iter, $\#(A - \sigma I)Y$, and $\#$eig\_out, refer to the caption of Table \ref{Tab:8-6-1}.
$s$ is the exact number of eigenvalues
inside an interval.
}
\footnotesize{
\noindent
\begin{tabular}{|c|cccccc|} \hline
Interval
& $s$ & $\#$eig\_out &$\tau_r$& $\tau_\lambda$ &$\#$iter &$\#(A - \sigma I)Y$\\ \hline
$(11.5, 12)$ &200&15 &2.07e-03 &1.07e-03 & 974&44 \\ \hline
$(11.6, 12)$ & 167&15 &1.93e-03 &2.20e-04 &960 &42 \\ \hline
$(11.7, 12)$ & 125&15 &3.47e-05& 6.08e-08 & 960&35\\ \hline
$(11.8, 12)$ &84 &15 &1.85e-06 &2.68e-10 &960 &31 \\ \hline
$(11.9, 12)$& 46 &15 &1.75e-06 &2.67e-10 & 960&3 \\ \hline
$(11.95, 12)$& 22&15 &3.00e-03& 1.52e-04&960&1 \\ \hline
$(11.99, 12)$& 5&5&2.34e-06&1.66e-10&960&0\\ \hline
$(11.995, 12)$&3& 3& 6.02e-07&4.37e-11&960&0  \\ \hline
$(11.998, 12)$ &1& 1&1.99e-06&1.54e-11&960&0  \\ \hline
$(11.999, 12)$ &0 &0 &- &- &960 &0 \\ \hline
\end{tabular}}
\label{Tab:9-11-1}
\end{table}

\begin{table}
\centering
\caption{Experiment 2: Results about Algorithm \ref{alg:7-29-1}. We set
$m = 60$, $num\_cmp = 30$, $num\_out = 15$, and $max\_it = 100$.
The meanings of $s$, $\#$eig\_out, $\tau_r$, $\tau_\lambda$, $\#$iter, and $\#(A - \sigma I)Y$ are
the same as
in the caption of Table \ref{Tab:9-11-1}.
}
\footnotesize{
\noindent
\begin{tabular}{|c|cccccc|} \hline
Interval
& $s$ & $\#$eig\_out &$\tau_r$& $\tau_\lambda$ &$\#$iter &$\#(A - \sigma I)Y$ \\ \hline
$(11.5, 12)$&200 &15 &1.59e-03 &1.28e-03 & 974&86  \\ \hline
$(11.6, 12)$&167 &15 & 6.08e-05& 2.86e-07&960 &82  \\ \hline
$(11.7, 12)$&125 &15  &4.32e-09  &3.86e-15 & 960&75 \\ \hline
$(11.8, 12)$&84 & 15 &6.76e-11  &3.56e-15 &960 &66 \\ \hline
$(11.9, 12)$&46 &15  &2.01e-11  & 3.11e-15 &960 &12  \\ \hline
$(11.95, 12)$& 22& 15 & 4.44e-09 & 3.26e-15 &960 &1  \\ \hline
$(11.99, 12)$&5 & 5 &2.77e-10  & 2.22e-15 & 960& 0\\ \hline
$(11.995, 12)$& 3&3&6.13e-11  & 4.00e-15 &960 &0  \\ \hline
$(11.998, 12)$&1 &1 &4.09e-11  &2.81e-15  &960 &0 \\ \hline
$(11.999, 12)$& 0& 0&-  &-  &960 &0 \\ \hline
\end{tabular}}
\label{Tab:9-14-1}
\end{table}

\begin{figure}
\begin{center}
\includegraphics[width=6.4cm]{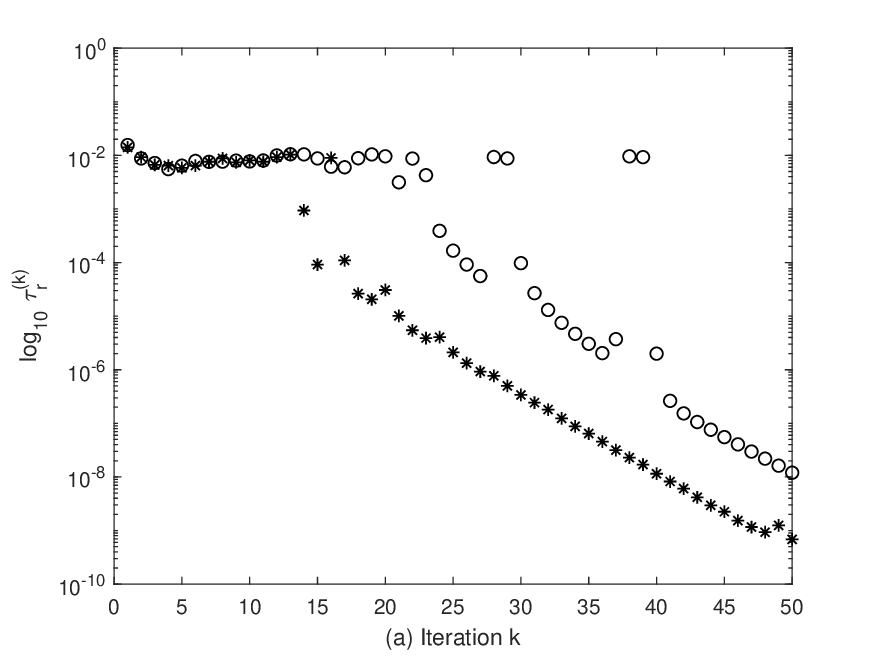}
\includegraphics[width=6.4cm]{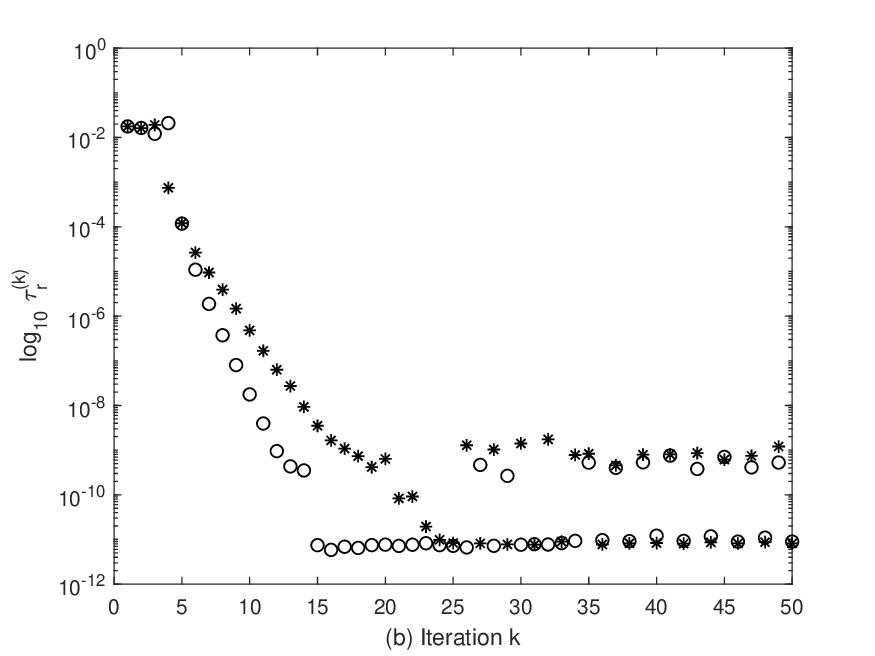}
\caption{Experiment 3: convergence histories of Algorithm \ref{alg:7-29-1}. $m = 60$, $num\_cmp = 30$, $num\_out = 30$.
(a) $(4.5, 5)$: $\circ$; $(4.6, 5)$: $*$.
(b) $(19.5, 20)$: $\circ$; $(19.6, 20)$: $*$.}
\label{Fig:10-10-1}
\end{center}
\end{figure}

{\bf Experiment 3.} In this experiment, we test the performance of Algorithm \ref{alg:7-29-1} at the two ends
of the spectrum. We select two intervals near each end, and fix $m = 60$, $num\_cmp = 30$, $num\_out = 30$ and $max\_it = 50$.
The experimental results are shown in Table \ref{Tab:10-9-1} and the histories of the relative residual norms $\tau_r^{(k)}$
 against iteration number $k$ are plotted in Figure \ref{Fig:10-10-1}.
The results reveal that Algorithm \ref{alg:7-29-1} converges faster near the right end of the spectrum. Probably it is because the spectrum
has a lower eigenvalue density at its right end,
resulting in relatively smaller ratios $|(\lambda_{m+1} - \sigma)/(\lambda_i - \sigma)|$ on which the convergence rate depends (see \S\ref{subsec: FFS}).\\

\begin{table}
\centering
\caption{Experiment 3: Results about Algorithm \ref{alg:7-29-1}. We set $m = 60$, $num\_cmp = 30$, $num\_out = 30$ and
$max\_it = 50$.
For the meanings of the quantities $s$, $\tau_r$, $\tau_\lambda$,
$\#$iter, $\#$eig\_out, and $\#(A - \sigma I)Y$, refer to the captions of Tables \ref{Tab:8-6-1} and \ref{Tab:9-11-1}. In each case of this experiment, $\#$eig\_out $= num\_out$.
}
\footnotesize{
\noindent
\begin{tabular}{|c|ccccc|} \hline
Interval & s & $\tau_r$& $\tau_\lambda$ &$\#$iter &$\#(A - \sigma I)Y$ \\ \hline
$(4.5, 5)$ & 118 &1.20e-08 &2.85e-14 &751 &53 \\ \hline
$(4.6, 5)$ & 90 &6.84e-10 & 4.21e-15& 751 &36 \\ \hline
$(19.5, 20)$ &49 &5.84e-12& 3.59e-15&947  &25 \\ \hline
$(19.6, 20)$ &40 &7.65e-12 &4.48e-15 &947 &25 \\ \hline
\end{tabular}}
\label{Tab:10-9-1}
\end{table}

{\bf Experiment 4.} We illustrate the scenario described in \S\ref{subsec:alleigs} of finding all the eigenvalues
in a given interval by Algorithm \ref{alg:7-29-1}.

Let us consider the interval $(a, b) = (11.7,12)$ and find all the $125$ eigenvalues in it. The length $\delta$ of this
interval is $0.3$. We set $m = 80$, $num\_cmp = 40$, $num\_out = 20$ and $max\_it = 100$
for Algorithm \ref{alg:7-29-1}. We first apply the algorithm
to the interval $(a, b)$ to obtain the first $20$
largest eigenvalues in it: $\hat{\lambda}^{(1)}_1 \geq \ldots \geq \hat{\lambda}^{(1)}_{20}$ with $\hat{\lambda}^{(1)}_1 = 11.9984$ and $\hat{\lambda}^{(1)}_{20}
= 11.9542$. Then pick a $b_1$. To this end, we evenly divide $(a, b)$ into ten subintervals and find that the subinterval $(11.94, 11.97)$ contains $\hat{\lambda}^{(1)}_{20}$. We then set $b_1 = 11.97$, $a_1 = b_1 - \delta = 11.67$, and apply Algorithm \ref{alg:7-29-1} to the interval $(a_1, b_1)$ to obtain the first $20$
largest eigenvalues in it: $\hat{\lambda}^{(2)}_1 \geq \ldots \geq \hat{\lambda}^{(2)}_{20}$ with $\hat{\lambda}^{(2)}_1 = 11.9686$ and $\hat{\lambda}^{(2)}_{20}
= 11.9278$. To pick a $b_2$, evenly divide $(a_1, b_1)$ into ten subintervals. Since $\hat{\lambda}^{(2)}_{20}$ lies in the subinterval $(11.91, 11.94)$, we set
$b_2 = 11.94$, and $a_2 = b_2 - \delta = 11.64$ accordingly, then apply Algorithm \ref{alg:7-29-1} to the interval $(a_2, b_2)$ to get the first $20$
largest eigenvalues in it: $\hat{\lambda}^{(3)}_1 \geq \ldots \geq \hat{\lambda}^{(3)}_{20}$ with $\hat{\lambda}^{(3)}_1 = 11.9369$ and $\hat{\lambda}^{(3)}_{20}
= 11.8977$. Each application of Algorithm \ref{alg:7-29-1} determines some, but not greater than $20$, eigenvalues.
We repeat this process until all the eigenvalues in $(a, b)$ have been found. Details of the numerical results are shown in Table \ref{Tab:8-4-101}.

\begin{table}
\centering
\caption{Experiment 4: Results about Algorithm \ref{alg:7-29-1}. We set
$m = 80$, $num\_cmp = 40$, $num\_out = 20$, and max\_it $= 100$.
The meanings of $s$, $\tau_r$, $\tau_\lambda$, $\#$iter, and $\#(A - \sigma I)Y$ are
the same as
in the caption of Table \ref{Tab:9-11-1}.
$\hat{\lambda}_{1}$ and $\hat{\lambda}_{20}$ are respectively the largest and
the $20$th largest eigenvalues in an interval, computed by Algorithm \ref{alg:7-29-1}.
}
\footnotesize{
\noindent
\begin{tabular}{|cc|ccccccc|} \hline
No. &Interval
& $s$ & $\hat{\lambda}_{20}$ & $\hat{\lambda}_{1}$ & $\tau_r$& $\tau_\lambda$ &$\#$iter &$\#(A - \sigma I)Y$ \\ \hline
0&$(11.7, 12)$ &125& 11.9542 &11.9984 &1.50e-11 &2.97e-15 & 961&74\\ \hline
1&$(11.67, 11.97)$& 120& 11.9278&11.9686&1.91e-11&5.65e-15 &971&97\\ \hline
2&$(11.64, 11.94)$ &119 & 11.8977&11.9369&1.68e-11&4.47e-15&959&61\\ \hline
3&$(11.61, 11.91)$&118 & 11.8596&11.9052&2.09e-11&2.69e-15&961&87\\ \hline
4&$(11.58, 11.88)$&117 &11.8230&11.8775&1.94e-11&5.10e-15&971&81 \\ \hline
5&$(11.55, 11.85)$ & 117 & 11.7932&11.8482&2.08e-11&4.50e-15&973&48\\ \hline
6&$(11.52, 11.82)$& 121 &11.7718&11.8183&1.85e-11&4.83e-15&960&84\\ \hline
7&$(11.49, 11.79)$& 116 &11.7406&11.7873&1.48e-11&4.07e-15&958&56\\ \hline
8&$(11.46, 11.76)$&113 &11.7189& 11.7585&2.11e-11 &4.39e-15&961& 82\\ \hline
9&$(11.43, 11.73)$& 113 & 11.6802& 11.7272&2.86e-11 &7.74e-15&961 &57\\ \hline
\end{tabular}}
\label{Tab:8-4-101}
\end{table}

We also report as in \cite{TP13} the orthogonality properties of the distinct computed eigenvectors in Table \ref{Tab:8-2-1}.
In the first column of Table \ref{Tab:8-4-101}, we have numbered the intervals. We then define $ortho_{ij} \equiv \max |x_i^T x_j|$ where $x_i$ and $x_j$ are the computed eigenvectors associated with the $i$th and the $j$th intervals respectively. For example, $ortho_{22} = 4.02 \times 10^{-12}$ is the maximum mutual orthogonality value of
the distinct computed eigenvectors from the interval $(a_2, b_2) = (11.64, 11.94)$. The overall orthogonality is $\displaystyle{ortho_{all} \equiv \max_{i, j} ortho_{ij} =
1.66 \times 10^{-9}}$.

\begin{table}
\centering
\caption{Experiment 4: (Continuation of Table \ref{Tab:8-4-101}) We report the mutual orthogonality of the computed eigenvectors where
$ortho_{ij} \equiv \max |x_i^T x_j|$ with $x_i$ and $x_j$ being the computed eigenvectors from the $i$th and the $j$th intervals respectively.
The overall orthogonality is $\displaystyle{ortho_{all} \equiv \max_{ij} ortho_{ij} =
1.66 \times 10^{-9}}$.
}
\footnotesize{
\noindent
\begin{tabular}{|c|c|c|c|c|c|} \hline
$ortho_{ij}$& 0 & 1 & 2 & 3&4 \\ \hline
0 & 2.77e-12 & 8.34e-12&8.87e-10&1.86e-11 &8.53e-12\\ \hline
1& 8.34e-12 &6.06e-12 &1.58e-11&1.07e-10 & 1.07e-10\\ \hline
2& 8.87e-10 &1.58e-11 &4.02e-12 &2.36e-11&6.06e-13\\ \hline
3& 1.86e-11&1.07e-10 &2.36e-11&6.72e-12&1.28e-11\\ \hline
4 & 8.53e-12&1.07e-10&6.06e-13&1.28e-11&2.05e-11\\ \hline
5& 8.38e-16&9.70e-16&5.24e-10&1.61e-09&1.76e-11\\ \hline
6 & 8.89e-16&9.97e-16&5.24e-10&1.69e-11&1.34e-12\\ \hline
7& 9.34e-16&8.77e-16&6.22e-13&1.47e-14&1.38e-10\\ \hline
8 & 6.73e-16&9.12e-16&1.05e-15&1.26e-13&4.87e-10\\ \hline
9 & 6.74e-16&8.14e-16&7.01e-16&8.76e-16&9.30e-16\\ \hline\hline
& 5 & 6 & 7 & 8&9 \\ \hline
0& 8.38e-16&8.89e-16 &9.34e-16&6.73e-16&6.74e-16\\ \hline
1& 1.00e-15 &9.97e-16 &8.42e-16 &9.48e-16 &8.14e-16\\ \hline
2&5.24e-10 &5.24e-10 &6.22e-13&1.05e-15&7.22e-16
\\ \hline
3&1.61e-09&1.69e-11 &1.47e-14& 1.26e-13&8.76e-16\\ \hline
4 &1.76e-11&1.34e-12 &1.38e-10&4.87e-10&9.71e-16\\ \hline
5& 7.28e-12&1.98e-11&1.38e-10&4.87e-10&1.08e-09\\ \hline
6 & 1.98e-11&8.35e-12&1.49e-11&4.34e-13&1.66e-09\\ \hline
7& 1.38e-10&1.49e-11&6.91e-12&1.61e-11&1.66e-09\\ \hline
8 & 4.87e-10&4.34e-13&1.61e-11&1.97e-11&3.60e-11\\ \hline
9 & 1.08e-09&1.66e-09&1.66e-09&3.60e-11&1.44e-11\\ \hline
\end{tabular}}
\label{Tab:8-2-1}
\end{table}

\subsection{Experiments with Andrews}\label{subsec:Andrews} For each of the three points $5, 18$ and $31$ on the real axis in the complex plane, we
use the \textsc{Matlab} command $[V, D] =$ \texttt{eigs}$(A, k, sigma)$ to find $500$ eigenvalues closest to the point, together
with their corresponding eigenvectors, of the {\it Andrews} matrix $A$.
The eigenvalues $\lambda_i$ and eigenvectors $v_i$ obtained
satisfy $\displaystyle{
\max_{1 \leq i \leq 500} \|Av_i - \lambda_i v_i\|_2/\|v_i\|_2 < 9.52 \times 10^{-10}, 4.06 \times 10^{-10}}$, and $1.16 \times 10^{-11}$ respectively. For this matrix, the scale factor $\mu \approx 12.59$.\\

{\bf Experiment 5.} We repeat Experiment 2 on the matrix {\it Andrews}. The sequence of intervals chosen and the detailed
numerical results are presented in Table \ref{Tab:11-11-1}. In this experiment, similar observations to those in Experiment 2 can be made.\\

\begin{table}
\centering
\caption{Experiment 5. Results about Algorithm \ref{alg:7-29-1}. We set $m = 80$, $num\_cmp = 40$, $num\_out = 20$ and $max\_it = 100$.
For the meanings of the quantities $s$, $\tau_r$, $\tau_\lambda$,
$\#$iter, $\#$eig\_out, and $\#(A - \sigma I)Y$, refer to the captions of Tables \ref{Tab:8-6-1} and \ref{Tab:9-11-1}.
}
\footnotesize{
\noindent
\begin{tabular}{|c|cccccc|} \hline
Interval  &s & $\tau_r$& $\tau_\lambda$ &$\#$iter &$\#$eig\_out &$\#(A - \sigma I)Y$ \\ \hline\hline
$(17.87, 18)$ &241 & 1.00e-04&5.90e-07 &3348 &20 &66 \\ \hline
$(17.91, 18)$ &164 &4.65e-08 &2.05e-13 & 3376& 20&68 \\ \hline
$(17.93, 18)$ &130 &1.31e-08 &3.20e-14 &3359 &20 &47 \\ \hline
$(17.95, 18)$ &91&7.83e-09 &1.50e-14  &3362 &20 &27 \\ \hline
$(17.97, 18)$ &53 &3.17e-07 &6.49e-12 &3359 &20 &4  \\ \hline
$(17.99, 18)$ &17 &9.11e-06 & 3.16e-09&3359 &17 &1 \\ \hline
$(17.995, 18)$ &10 &3.29e-07 & 2.46e-12&3369 &10 &0 \\ \hline
$(17.999, 18)$ &3 &6.44e-08 &6.05e-13 &3364 &3 &0 \\ \hline
$(17.9995, 18)$ &1 &2.85e-07 & 2.09e-12&3382 & 1& 0\\ \hline
$(17.9999, 18)$ &0 &- &- &3376 &0 & 0\\ \hline
\end{tabular}}
\label{Tab:11-11-1}
\end{table}

\begin{figure}
\begin{center}
\includegraphics[width=6.4cm]{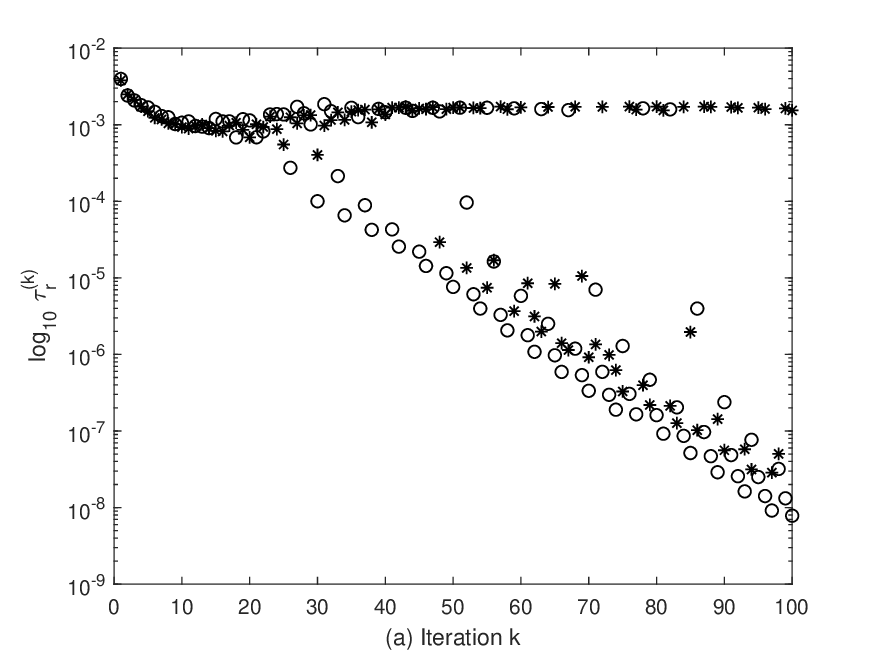}
\includegraphics[width=6.4cm]{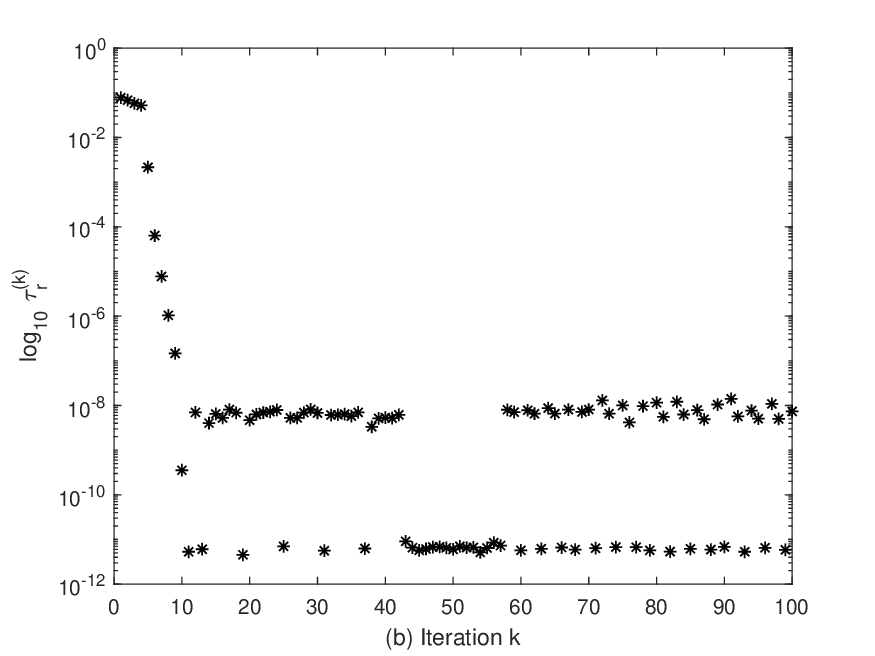}
\caption{Experiment 6: convergence histories of Algorithm \ref{alg:7-29-1}.
(a) $(4.95, 5)$: $*$; $(17.95, 18)$: $o$. (b) $(30, 31)$: $*$.}
\label{Fig:10-4-1}
\end{center}
\end{figure}

{\bf Experiment 6.}
We test Algorithm \ref{alg:7-29-1} on the three intervals $(4.95, 5), (17.95, 18)$, and $(30, 31)$, locating at the two ends and in the middle of the spectrum of $A$ respectively. The numerical results are shown in Table \ref{Tab:10-11-1}, and the convergence histories of $\tau_r^{(k)}$ against iteration number $k$ are plotted in Figure \ref{Fig:10-4-1}.

Among the three intervals, eigenpairs in $(4.95, 5)$ are the most difficult to compute.
The relative residual $\tau_r^{(k)}$ remains about
$O(10^{-3})$
in the first $45$ iterations before it
starts to drop
(see Figure \ref{Fig:10-4-1}(a)).\\

\begin{table}
\centering
\caption{Experiment 6: Results about Algorithm \ref{alg:7-29-1}.
We set $max\_it =
100$.
For the meanings of the quantities $m$, nc, no, $s$, $\tau_r$, $\tau_\lambda$,
$\#$iter, $\#$eig\_out, and $\#(A - \sigma I)Y$, refer to the captions of Tables \ref{Tab:8-6-1} and \ref{Tab:9-11-1}. In each case of this experiment, $\#$eig\_out $= num\_out$.
}
\footnotesize{
\noindent
\begin{tabular}{|cccc|ccccc|} \hline
Interval & m & nc & no &s & $\tau_r$& $\tau_\lambda$ &$\#$iter &$\#(A - \sigma I)Y$ \\ \hline
$(4.95, 5)$ &100 &50 &25 & 113& 2.86e-08& 4.42e-13 &2318 &28 \\ \hline
$(17.95, 18)$ & 80 &40 & 20  & 91&7.83e-09 &1.50e-14  &3362 &27 \\ \hline
$(30, 31)$ &50 &25 & 25&26 &4.54e-12 & 1.58e-13&2123 &64 \\ \hline
\end{tabular}}
\label{Tab:10-11-1}
\end{table}

\subsection{Experiment with a random matrix}\label{sec:10-2-1} The experiment is motivated by the fact that a symmetric $A \in {\mathbb R}^{n \times n}$
is orthogonally diagonalizable.

Let $n = 10^6$ and $\eta = n/1500$. We generate a random diagonal matrix $A$ in \textsc{Matlab} with independent and uniformly distributed diagonal entries from the interval
$[-0.5\eta, 0.5\eta]$.
The ANE of such a matrix
is $n/\eta = 1500$, about the same as the ANE of the {\it Andrews} matrix. The scale factor $\mu$ in (\ref{equ:7-30-1}) is about $192.59$.

BiCG can solve the linear systems in (\ref{equ:7-15-3}) with a maximum number of iterations being about $65000$.
However, it took too long to complete this sequential process.
Because of the diagonal structure of $A$, we then decided to solve the linear systems with the \textsc{Matlab} operator ``$./$'' instead of
using BiCG. Specifically, consider the diagonal linear system
\begin{equation}\label{equ:1-4-1}
[(c + r e^{i \pi t_k}) I - A] x = y_j.
\end{equation}
Let $u$ be the vector of the main diagonal entries of the coefficient matrix in (\ref{equ:1-4-1}), and denote the $i$th entries of the vectors $u$ and $y_j$ by $u_i$ and $y_{ji}$ respectively.
Then we compute $x = [x_1, x_2, \ldots, x_n]^T = [y_{j1} / u_1, y_{j2} / u_2, \ldots, y_{jn} / u_n]^T$. After $x$ is computed, we
add some small perturbation to $x$ to mimic the solution of the system
(\ref{equ:1-4-1}) by BiCG:
\begin{equation}\label{equ:12-19-1}
x = x + \xi * \texttt{randn}(n, 1)
\end{equation}
where $\xi = 10^{-10} \|y_j\|_2/\|u\|_2$. It can be seen that the perturbed $x$ in (\ref{equ:12-19-1}) satisfies
$$
\|[(c + r e^{i \pi t_k}) I - A] x - y_j\|_2/\|y_j\|_2 \approx 10^{-10}.
$$
In this experiment, we use the $x$ in (\ref{equ:12-19-1}) as the numerical solution to the system (\ref{equ:1-4-1}).\\

{\bf Experiment 7.} We repeat Experiment 2 on the random matrix $A$. The
intervals chosen and the detailed
numerical results are presented in Table \ref{Tab:12-20-5}.
Besides the observations similar to those in Experiment 2, we also note that
a reasonable value for $r$ depends on the ANE of a matrix. In this experiment,
we choose $r = 2$
which yields
Table \ref{Tab:12-20-5}.
When we chose $r = 3$, however, Algorithm \ref{alg:7-29-1} did not converge well within
$100$ iterations.

\begin{table}
\centering
\caption{Experiment 7:
Results about Algorithm \ref{alg:7-29-1}.
We set $m = 80$, $num\_cmp = 40$, $num\_out = 20$ and $max\_it = 100$.
For the meanings of the quantities $s$, $\tau_r$, $\tau_\lambda$,
$\#$eig\_out, and $\#(A - \sigma I)Y$, refer to the captions of Tables \ref{Tab:8-6-1} and \ref{Tab:9-11-1}.
}
\footnotesize{
\noindent
\begin{tabular}{|c|ccccc|} \hline
Interval  & s&$\tau_r$& $\tau_\lambda$ &$\#eig\_out$ &$\#(A - \sigma I)Y$\\ \hline
(1, 1.21)&307 &3.34e-05 &4.22e-03 &20 &97  \\ \hline
(1, 1.19)&277 &3.53e-05 &7.30e-03 &20 & 74  \\ \hline
(1, 1.17)&249&4.65e-08 &1.64e-09 &20 &98  \\ \hline
(1, 1.15)&212 & 8.74e-08&1.39e-09 &20 &86  \\ \hline
(1, 1.13)&187 & 3.58e-08& 5.48e-10&20 &68 \\ \hline
(1, 1.11)& 155&3.60e-11 &4.62e-15 &20 &78 \\ \hline
(1, 1.09)& 125&8.02e-11 &6.15e-15 &20& 62\\ \hline
(1, 1.07)& 97&9.14e-13 &5.42e-15 & 20&43\\ \hline
(1, 1.05)& 72& 3.37e-13&5.74e-15&20&21\\ \hline
(1, 1.03)&38&9.88e-09&7.32e-12&20&1\\ \hline
(1, 1.01)& 18&6.25e-12&4.42e-15 &18&0 \\ \hline
(1, 1.001)& 2& 8.64e-10  &3.98e-13 &2&1 \\ \hline
(1, 1.0002)& 1&1.45e-12  &2.22e-16 &1&0 \\ \hline
(1, 1.0001)& 0& - &- & 0&0\\ \hline
\end{tabular}}
\label{Tab:12-20-5}
\end{table}


\section{Conclusions}\label{sec:conclusions} We incorporate a power subspace iteration process into the FEAST eigensolver to solve real and symmetric eigenvalue problems.
Together with two contour integrations per iteration, our approach has the advantages described at the end of \S\ref{sec:intro}.
Numerical experiments show that the resulting algorithm $F_2P$
is a robust and
accurate eigensolver for the computation of extreme as well as interior eigenvalues. We also observe that
$F_2P$ does not require $m \ge s$, but
$m$ should be chosen near
$s$.
 Moreover,
 it seems that there is some relation between the circle radius $r$ and the average number of eigenvalues per unit interval.

More experiments, especially on test data of large size
(e.g., hundred thousands
or more),
are needed to better understand the behavior of the
algorithm.
%
Our future work includes further reducing the dependence of $m$ on $s$ and applying $F_2P$ to the solution
of extremely ill-conditioned linear systems.




\begin{thebibliography}{10}







\bibitem{bks}
{\sc C. Bekas, E. Kokiopoulou, and Y. Saad}, {\it Computation of large invariant subspaces using polynomial filtered Lanczos
iterations with applications in density functional theory}, SIAM J. Matrix Anal. and Appl., 1(2008), pp. 397-418.









\bibitem{ufl_sparse_matrics}
{\sc T. A. Davis and Y. Hu}, {\it The University of Florida Sparse Matrix Collection},
ACM Transactions on Mathematical Software, 38:1--25, 2011.
http://www.cise.ufl.edu/research/sparse/matrices.

\bibitem{DR84}
{\sc P.~J. Davis and P. Rabinowitz}, {\it Methods of numerical integration}, 2nd Edition, Academic Press, Orlando, FL, 1984.



\bibitem{FangSaad}
{\sc H. R. Fang and Y. Saad}, {\it A filtered Lanczos procedure for extreme and interior eigenvalue problems},
SIAM J. Sci. Comput., 34 (2012), pp. A2220-A2246.

\bibitem{Fletcher}
{\sc R. Fletcher}, {\it Conjugate gradient methods for indefinite systems}, In volume 506 of Lecture Notes
Math., pages 73-89. Springer-Verlag, Berlin-Heidelberg-New York, 1976.


\bibitem{futa}
{\sc Y. Futamura, H. Tadano, and T. Sakurai}, {\it Parallel stochastic estimation method of eigenvalue distribution}, JSIAM Letters 2 (2010), pp.127--130.



\bibitem{GP16}
{\sc B. Gavin and E. Polizzi}, {\it Enhancing the performance and robustness
of the feast eigensolver}, In High Performance Extreme Computing Conference
(HPEC), 2016 IEEE (2016), IEEE, pp. 1-6.

\bibitem{GP}
{\sc ------},
{\it Krylov eigenvalue strategy using the FEAST algorithm with
inexact system solves}, Numer. Linear Algebra Appl. 25(2018), no. 5, e2188.

\bibitem{GPTV}
{\sc S. G\"{u}ttel, E. Polizzi, P. T. P. Tang, and G. Viaud}, {\it Zolotarev quadrature rules and load balancing
for the FEAST eigensolver}, SIAM J. Sci. Comput., 37(2015), pp. A2100-A2122.

\bibitem{gvl}
{\sc G. H. Golub and C. F. Van Loan}, {\it Matrix Computations}, 3rd Edition, Johns Hopkins University Press, Baltimore, MD, 1996.








\bibitem{kpt}
{\sc J. Kestyn, E. Polizzi, and P. T. P. Tang}, {\it FEAST eigensolver for non-Hermitian problems}, SIAM J. Sci. Comput., 38(5):S772-S799, 2016.

\bibitem{kbn}
{\sc K. Kollnig, P. Bientinesi, and E. D. Napoli},
{\it Rational spectral filters with optimal convergence rate}, SIAM J. Sci. Comput., 43(4), A2660-A2684, 2021.

\bibitem{kramer}
{\sc L. Kr\"amer, E. D. Napoli, M. Galgon, B. Lang, and P. Bientinesi}, {\it Dissecting the FEAST algorithm for generalized eigenproblems}, J. Comput. Appl. Math., 244 (2013), pp. 1--9.






\bibitem{NPY}
{\sc E. D. Napoli, E. Polizzi, and Y. Saad}, {\it Efficient estimation of eigenvalue counts in an interval}, Numer. Linear Algebra Appl. 2016; 23:674-692.

\bibitem{ps}
{\sc C. C. Paige and M. A. Saunders}, {\it Solution of sparse indefinite systems of linear equations}, SIAM J.
Numer. Anal., 12 (1975), pp. 617-629.

\bibitem{Parlett}
{\sc B. N. Parlett}, {\it The Symmetric Eigenvalue Problem}, no. 20 in Classics in Applied Mathematics, SIAM Publications,
Philadelphia, PA, 1998.

\bibitem{polizzi}
{\sc E. Polizzi}, {\it Density-matrix-based algorithm for solving eigenvalue problems}, Phys. Rev. B 79 (2009) 115112.

\bibitem{polizzi_software}
{\sc E. Polizzi et al.}, {\it FEAST eigenvalue solver}, http://www.feast-solver.org/.



\bibitem{saad03}
{\sc Y. Saad}, {\it Iterative Methods for Sparse Linear Systems}, SIAM, Philadelphia, 2nd edition, 2003.

\bibitem{saad}
{\sc ------},
{\it Numerical Methods for Large Eigenvalue Problems}, SIAM, Philadelphia, 2011.

\bibitem{gmres}
{\sc Y. Saad and M. H. Schultz}, {\it GMRES: A generalized minimal residual algorithm for solving nonsymmetric linear systems}, SIAM J. Sci. Stat. Comput., 7:856-869, 1986.





 \bibitem{SGKRH}
 {\sc O. Schenk, K. G\"artner, G. Karypis, S. R\"ollin, and M. Hagemann}, {\it PARDISO solver project}, 2010. http://www.pardiso-project.org/





\bibitem{TP13}
{\sc P.~T.~P. Tang and E. Polizzi}, {\it FEAST as a subspace iteration eigensolver accelerated by approximate spectral projection}, SIAM J. Matrix Anal. Appl.,  35 (2014), pp. 354--390.


\bibitem{Vorst}
{\sc H. A. Van der Vorst}, {\it Iterative Krylov methods for large linear systems},
Cambridge University Press,
2003.

\bibitem{Viaud}
{\sc G. Viaud}, {\it The FEAST Method}, M. Sc. dissertation, University of Oxford, 2021.


\bibitem{bib:WD2019}
{\sc J.~Winkelmann and E.~Di Napoli},
{\it Non-linear least-squares optimization of rational filters for the solution of interior Hermitian eigenvalue problems},
Frontiers in Applied Mathematics and Statistics, 5(5) (2019), pp. 1--17.


\bibitem{xisaad}
{\sc Y. Xi and Y. Saad},
{\it Computing partial spectra with least-squares rational
filters}, SIAM J Sci Comput. (2016) 38:A3020-45.

\bibitem{yxccb}
{\sc X. Ye, J. Xia, R. Chan, S. Cauley, and V. Balakrishnan}, {\it A Fast Contour-Integral Eigensolver for Non-Hermitian Matrices}, SIAM J. Matrix Anal. Appl., 38, (2017), pp. 1268-1297.

\bibitem{yin2019}
{\sc G.~J. Yin}, {\it A harmonic FEAST algorithm for non-Hermitian generalized eigenvalue problems},  Linear Algebra Appl., 578 (2019), pp. 75--94.

\bibitem{ycy}
{\sc G. Yin, R. Chan, and M. Yeung},  {\it A FEAST algorithm with oblique projection for generalized eigenvalue problems},
Numerical Linear Algebra with Applications, 2017; 24:e2092.
\end{thebibliography}
\end{document}